\documentclass{article}

\usepackage[utf8]{inputenc} 
\usepackage[T1]{fontenc}    
\usepackage{url}            
\usepackage{booktabs}       
\usepackage{amsfonts}       
\usepackage{microtype}      
\usepackage{nicefrac}
\usepackage{bbm}

\usepackage{algpseudocode}
\usepackage{algorithm}

\usepackage{bm}
\usepackage[colorinlistoftodos]{todonotes}
\usepackage{amsmath}
\usepackage{amsthm}
\usepackage{amssymb}
\usepackage{thmtools}
\usepackage{thm-restate}
\usepackage{paralist}
\usepackage{comment}
\usepackage{cite}
\usepackage{float}
\usepackage{balance}
\usepackage{stmaryrd}
\usepackage{subcaption}
\usepackage{mathtools}

\usepackage{enumitem}
\usepackage{fancyhdr}
\usepackage{wrapfig}

\usepackage{multirow}
\usepackage{pifont}

\usepackage[colorlinks,urlcolor=blue,citecolor=blue,linkcolor=blue]{hyperref} 

\newcommand{\innp}[1]{\left\langle #1 \right\rangle}

\newcommand{\vx}{\mathbf{x}}
\newcommand{\vd}{\mathbf{d}}

\newcommand{\cx}{\mathcal{X}}
\newcommand{\cc}{\mathcal{C}}
\newcommand{\cs}{\mathcal{S}}
\newcommand{\cl}{\mathcal{L}}

\newcommand{\vy}{\mathbf{y}}
\newcommand{\vz}{\mathbf{z}}

\newcommand{\vw}{\mathbf{w}}
\newcommand{\vr}{\mathbf{r}}

\newcommand{\va}{\mathbf{a}}

\newcommand{\vu}{\mathbf{u}}
\newcommand{\vs}{\mathbf{s}}

\newcommand{\vlambda}{\bm{\lambda}}

\newcommand{\defeq}{\stackrel{\mathrm{\scriptscriptstyle def}}{=}}

\newcommand{\rr}{\mathbb{R}}
\newcommand{\norm}[1]{\left\| #1 \right\|}

\newcommand{\boundary}[1]{\operatorname{Bd}(#1)}

\newcommand*{\vsepfbox}[1]{%
  \begingroup
    \sbox0{\fbox{#1}}%
    \setlength{\fboxrule}{0pt}%
    \mbox{\kern-\fboxsep\fbox{\unhbox0}\kern-\fboxsep}%
  \endgroup
}

\newcommand{\cmark}{\ding{51}}%
\newcommand{\xmark}{\ding{55}}%

\theoremstyle{plain} \numberwithin{equation}{section}
\newtheorem{theorem}{Theorem}[section]
\numberwithin{theorem}{section}

\newtheorem{lemma}[theorem]{Lemma}
\newtheorem{proposition}[theorem]{Proposition}

\theoremstyle{definition}
\newtheorem{definition}[theorem]{Definition}

\newtheorem{remark}[theorem]{Remark}
\newtheorem{example}[theorem]{Example}

\DeclareMathOperator{\interior}{\mathrm{Int}}

\DeclareMathOperator*{\argmin}{argmin}
\DeclareMathOperator*{\argmax}{argmax}

\newcommand{\dom}[1]{
  \mathrm{dom}(#1)
}

\graphicspath{{imgs/}}
\makeatletter
\def\mathcolor#1#{\@mathcolor{#1}}
\def\@mathcolor#1#2#3{%
  \protect\leavevmode
  \begingroup
    \color#1{#2}#3%
  \endgroup
}
\makeatother

\newcommand{\hrulealg}[0]{\vspace{1mm} \hrule \vspace{1mm}}

\usepackage[square,comma, sort&compress]{natbib}
\usepackage[margin=1in]{geometry}
\usepackage{jmlr2e}

\usepackage{charter}
\usepackage[vvarbb]{newtxmath}
\usepackage[capitalize,nameinlink]{cleveref}[0.21]

\makeatletter
\newenvironment{proof*}[1][\proofname]{\par
  \pushQED{\qed}%
  \normalfont \partopsep=\z@skip \topsep=\z@skip
  \trivlist
  \item[\hskip\labelsep
        \itshape
    #1\@addpunct{.}]\ignorespaces
}{%
  \popQED\endtrivlist\@endpefalse
}
\makeatother

\title{Scalable Frank-Wolfe on Generalized Self-concordant Functions via Simple Steps}

\author{\name Alejandro Carderera \email  \href{mailto:alejandro.carderera@gatech.edu}{alejandro.carderera@gatech.edu}\\
       \addr Department of Industrial and Systems Engineering\\
       Georgia Institute of Technology\\
       Atlanta, USA
       \AND
       \name Mathieu Besançon \email \href{mailto:mathieu.besancon@polymtl.ca}{besancon@zib.de} \\
       \addr Zuse Institute Berlin, Germany, and Univ.~Grenoble Alpes, Inria, LIG, Grenoble, France \\
       \AND
       \name Sebastian Pokutta \email \href{mailto:pokutta@zib.de}{pokutta@zib.de} \\
       \addr Institute of Mathematics\\
       Technische Universität Berlin and Zuse Institute Berlin, Germany}

\begin{document}

\maketitle

\begin{abstract} 
Generalized self-concordance is a key property present
in the objective function of many important learning problems.  We establish
the convergence rate of a simple Frank-Wolfe variant that uses the open-loop
step size strategy $\gamma_t = 2/(t+2)$, obtaining a $\mathcal{O}(1/t)$
convergence rate for this class of functions in terms of primal gap and
Frank-Wolfe gap, where $t$ is the iteration count. This avoids the use of
second-order information or the need to estimate local smoothness parameters of
previous work. We also show improved convergence rates for various common
cases, e.g., when the feasible region under consideration is uniformly convex
or polyhedral.  
\end{abstract}

\section{Introduction}

Constrained convex optimization is the cornerstone of many machine learning problems.
We consider such problems, formulated as:
\begin{equation}
\label{prob:main}
\min_{\vx \in \cx} f(\vx),
\end{equation}
where $f: \mathbb{R}^n \rightarrow \mathbb{R} \cup \{+\infty\}$ is a generalized self-concordant function and $\cx \subseteq \rr^n$ is a compact convex set.
When computing projections  onto the feasible regions as required in, e.g., projected gradient descent, is prohibitive, Frank-Wolfe (FW) \citep{frank1956algorithm} algorithms (a.k.a.~Conditional Gradients (CG) \citep{polyak66cg}) are the algorithms of choice, relying on Linear Minimization Oracles (LMO) at each iteration to solve Problem~\eqref{prob:main}. The analysis of their convergence often relies on the assumption that the gradient is Lipschitz-continuous. 
This assumption does not necessarily hold for generalized self-concordant functions, an important class of functions whose growth can be unbounded.
\subsection{Related work}

In the classical analysis of Newton's method, when the Hessian of $f$ is assumed to be Lipschitz continuous and the function is strongly convex, one arrives at a convergence rate for the algorithm that depends on the Euclidean structure of $\rr^n$, despite the  fact that the algorithm is affine-invariant. This motivated the introduction of self-concordant functions in \cite{nesterov1994interior}, functions for which the third derivative is bounded by the second-order derivative, with which one can obtain an affine-invariant convergence rate for the aforementioned algorithm.
More importantly, many of the barrier functions used in interior-point methods are self-concordant, which extends the use of polynomial-time interior-point methods to many settings of interest.

Self-concordant functions have received strong interest in recent years due to the attractive properties that they allow to prove for many statistical estimation settings \citep{pmlr-v99-marteau-ferey19a,ostrovskii2021finite}. The original definition of self-concordance has been expanded and generalized since its inception, as many objective functions of interest have self-concordant-like properties without satisfying the strict definition of self-concordance. For example, the logistic loss function used in logistic regression is not strictly self-concordant, but it fits into a class of pseudo-self-concordant functions, which allows one to obtain similar properties and bounds as those obtained for self-concordant functions \citep{bach2010self}. This was also the case in \citet{ostrovskii2021finite} and \cite{tran2015composite2}, in which more general properties of these
pseudo-self-concordant functions were established. This was fully formalized in \cite{sun_generalized_2019}, in which the concept of \emph{generalized self-concordant} functions was introduced, along with key bounds, properties, and variants of Newton methods for the unconstrained setting which make use of this property.

Most algorithms that aim to solve Problem~\eqref{prob:main} assume access to
second-order information, as this often allows the algorithms to make
monotonic progress, remain inside the domain of $f$, and often, converge quadratically when close enough to the optimum.
Recently, several lines of work have focused on using Frank-Wolfe algorithm variants to solve these types of problems in the projection-free
setting, for example constructing second-order approximations to a
self-concordant $f$ using first and second-order information, and minimizing
these approximations over $\cx$ using the Frank-Wolfe algorithm
\citep{liu2020newton}. Other approaches, such as the ones presented in
\citet{dvurechensky2020self} (later extended in
\citet{dvurechensky2020generalized}), apply the Frank-Wolfe algorithm to a
generalized self-concordant $f$, using first and second-order information about
the function to guarantee that the step sizes are so that the iterates do not
leave the domain of $f$, and monotonic progress is made. An additional Frank-Wolfe
variant in that work, in the spirit of
\citet{garber2016linearly}, utilizes first and second
order information about $f$, along with a Local Linear Optimization Oracle
for $\cx$, to obtain a linear convergence rate in primal gap over polytopes given in inequality description. The authors in \cite{dvurechensky2020generalized} also present an 
additional Frank-Wolfe variant which does not use second-order information, and 
uses the backtracking line search of
\citet{pedregosa2018step} to estimate local smoothness parameters at a given
iterate. Other specialized Frank-Wolfe algorithms
have been developed for specific problems involving generalized self-concordant
functions, such as the Frank-Wolfe variant developed for marginal inference with concave
maximization \citep{krishnan2015barrier}, the variant developed
in \cite{zhao_analysis_2020} for $\theta$-homogeneous barrier functions,
or the application for phase retrieval in \cite{odor2016frank}, where the Frank-Wolfe algorithm is shown to converge
on a self-concordant non-Lipschitz smooth objective.

\subsection{Contribution}

The contributions of this paper are detailed below and summarized in Table~\ref{tab:results}.
\paragraph{Simple FW variant for generalized self-concordant functions} We show that a small variation of the original Frank-Wolfe algorithm \citep{frank1956algorithm} with an open-loop step size of the form $\gamma_t = 2/(t+2)$, where $t$ is the iteration count is all that is needed to achieve
a convergence rate of $\mathcal{O}(1/t)$ in primal gap; this also answers an open question posed in \cite{dvurechensky2020generalized}. Our variation ensures monotonic progress while employing an open-loop strategy which, together with the iterates being convex combinations, ensures that we do not leave the domain of $f$. In contrast to other methods that depend on either a line search or second-order information, our variant uses only a linear minimization oracle, zeroth-order and first-order information and a domain oracle for $f(\vx)$. The assumption of the latter oracle is very mild and was also implicitly assumed in several of the algorithms presented in \cite{dvurechensky2020generalized}.
As such, our iterations are much cheaper than those in previous work, while essentially achieving the same convergence rates for Problem~\eqref{prob:main}.

Moreover, our variant relying on the open-loop step size $\gamma_t = 2/(t+2)$ allows us to establish a $\mathcal{O}(1/t)$ convergence rate  for the Frank-Wolfe gap, is agnostic, i.e., does not need to estimate local smoothness parameters, and is parameter-free, leading to convergence rates and oracle complexities that are independent of any tuning parameters.

\begin{table}[H]\centering
  \footnotesize
\begin{tabular}{lccccc}
\toprule
\multirow{2}{*}{\textbf{Algorithm}} & \multicolumn{2}{c}{\textbf{Convergence}} & \multirow{2}{*}{\textbf{Reference}} & \textbf{$1^{\text{st}}$-order} / & \multirow{2}{*}{\textbf{Requirements}}  \\
				    & \textbf{Primal gap}   & \textbf{FW gap}   & & \textbf{LS free?} & \\
\midrule
\texttt{FW-GSC}  & $\mathcal{O}(1/\varepsilon)$ &  & \cite[Alg.2]{dvurechensky2020generalized} & {\color{red} \xmark} /  {\color{green} \cmark}  &  SOO \\
\texttt{LBTFW-GSC}  & $\mathcal{O}(1/\varepsilon)$ &  & \cite[Alg.3]{dvurechensky2020generalized} & {\color{green} \cmark} /  {\color{red} \xmark}  & ZOO, DO  \\
\texttt{MBTFW-GSC}  & $\mathcal{O}(1/\varepsilon)$ &  & \cite[Alg.5]{dvurechensky2020generalized} & {\color{red} \xmark} /  {\color{green} \cmark}  &  ZOO, SOO, DO \\
\texttt{FW-LLOO} & $\mathcal{O}(\log 1/\varepsilon)$ &  & \cite[Alg.7]{dvurechensky2020generalized} & {\color{red} \xmark} /  {\color{green} \cmark} & $P(\mathcal{X})$, LLOO, SOO \\
\texttt{ASFW-GSC}  & $\mathcal{O}(\log 1/\varepsilon)$ &  & \cite[Alg.8]{dvurechensky2020generalized} & {\color{red} \xmark} /  {\color{green} \cmark}  &  $P(\mathcal{X})$, SOO \\
\texttt{M-FW} & $\mathcal{O}(1/\varepsilon)$ & $\mathcal{O}(1/\varepsilon)$ & \textbf{This work} & {\color{green} \cmark} /  {\color{green} \cmark} & ZOO, DO  \\
\texttt{B-\{AFW/BPCG\}} & $\mathcal{O}(\log 1/\varepsilon)$ & $\mathcal{O}(\log 1/\varepsilon)$ &  \textbf{This work} & {\color{green} \cmark}  /  {\color{red} \xmark} & $P(\mathcal{X})$, ZOO, DO \\
\bottomrule
\end{tabular}
\smallskip 

\caption{\label{tab:results}
Number of iterations needed to achieve an $\varepsilon$-optimal solution for Problem~\ref{prob:main}. We denote line search by LS, zeroth-order oracle by ZOO, second-order oracle by SOO, domain oracle by DO, local linear optimization oracle by LLOO, and the assumption that $\mathcal X$ is polyhedral by $P(\mathcal{X})$. The oracles listed under the \textbf{Requirements} column are the additional oracles required, other than the first-order oracle (FOO) and the linear minimization oracle (LMO) which all algorithms use.}
 \label{table:comparison}
\end{table}

\paragraph{Faster rates in common special cases} We also obtain improved
convergence rates when the optimum is contained in the interior of $\cx \cap \dom f$,
or when the set $\cx$ is uniformly or strongly convex, using the
backtracking line search of \cite{pedregosa2018step}.
We also show that the Away-step Frank-Wolfe \cite{wolfe70,lacoste15} and the Blended Pairwise Conditional Gradient \cite{tsuji2022pairwise} can use the aforementioned line search to achieve linear rates over polytopes.
For clarity we want to stress that any linear rate over polytopes has to depend also on the ambient dimension of the polytope; this applies to our linear rates and those in Table~\ref{table:comparison} established elsewhere (see \citet{CDP2019}).
In contrast, the $\mathcal O(1/\varepsilon)$ rates are dimension-independent.

\paragraph{Numerical experiments}
We provide numerical experiments that showcase the performance of the algorithms on generalized self-concordant objectives to complement the theoretical results.
In particular, they highlight that the simple step size strategy we propose is competitive with and sometimes outperforms other variants on many instances.

After publication of our initial draft, in a revision of their original work, \cite{dvurechensky2020generalized} added an analysis of the Away-step Frank-Wolfe algorithm which is complementary to ours (considering a slightly different setup and regimes) and was conducted independently; we have updated the tables to include these additional results.

\subsection{Preliminaries and Notation}

We denote the \emph{domain} of $f$ as $\dom f \defeq \{\vx \in \rr^n, f(\vx) < +\infty\}$ and the (potentially non-unique) minimizer of Problem~\eqref{prob:main} by $\vx^*$. Moreover, we denote the \emph{primal gap} and the \emph{Frank-Wolfe gap} at $\vx \in \cx \cap \dom f$ as $h(\vx) \defeq f(\vx) - f(\vx^*)$ and $g(\vx) \defeq \max_{\mathbf{v} \in \cx } \innp{\nabla f(\vx), \vx - \mathbf{v}}$, respectively. We use $\norm{\cdot}$, $\norm{\cdot}_{H}$, and $\innp{\cdot, \cdot}$ to denote the \emph{Euclidean norm}, the \emph{matrix norm} induced by a symmetric positive definite matrix $H\in \rr^{n\times n}$, and the \emph{Euclidean inner product}, respectively.
We denote the \emph{diameter} of $\cx$ as $D \defeq \max_{\vx, \vy \in \cx} \norm{\vx - \vy}$.
Given a non-empty set $\cx\subset \rr^n$  we refer to its \emph{boundary} as $\boundary{\cx}$ and to its \emph{interior} as $\interior\left(\cx\right)$.
We use $\Delta_n$ to denote the \emph{probability simplex} of dimension $n$. Given a compact convex set $\cc\subseteq \dom f$ we denote:
\begin{align*}
L_f^{\cc} = \max\limits_{\vu \in \cc, \vd \in \rr^n} \frac{\norm{\vd}_{\nabla^2 f(\vu)}^2}{\norm{\vd}^2_2},\;\;\;
\mu_f^{\cc} = \min\limits_{\vu \in \cc, \vd \in \rr^n} \frac{\norm{\vd}_{\nabla^2 f(\vu)}^2}{\norm{\vd}^2_2}.
\end{align*}
We assume access to:
\begin{enumerate}
\item \textbf{Domain Oracle (DO):} Given $\vx \in \cx$, return whether $\vx \in \dom f$.
\item \textbf{Zeroth-Order Oracle (ZOO)}: Given $\vx \in \dom f$, return $f(\vx)$.
\item \textbf{First-Order Oracle (FOO)}: Given $\vx \in \dom f$, return $\nabla f(\vx)$.
\item \textbf{Linear Minimization Oracle (LMO)}: Given $\vd \in \rr^n$, return $\argmin\limits_{\vx \in \cx}\innp{\vx, \vd}$.
\end{enumerate}
The FOO and LMO oracles are standard in the FW literature. The ZOO oracle is often implicitly assumed to be included with the FOO oracle; we make this explicit here for clarity. Finally, the DO oracle is motivated by the properties of generalized self-concordant functions. It is reasonable to assume the availability of the DO oracle: following the definition of the function codomain, one could simply evaluate $f$ at $\vx$ and assert $f(\vx) < +\infty$, thereby combining the DO and ZOO oracles into one oracle.
However, in many cases testing the membership of $\vx \in \dom f$ is computationally less demanding than the function evaluation.

\begin{remark}
	Requiring access to a zeroth-order and domain oracle are mild assumptions, that were also implicitly assumed in one of the three FW-variants presented in \citet{dvurechensky2020generalized} when computing the step size according to the strategy from \citet{pedregosa2018step}; see \cref{backtrack:domTestBacktrack} in \cref{backtrack}. The remaining two variants ensure that $\vx \in \dom f$ by using second-order information about $f$, which we explicitly do not rely on.
\end{remark}

The following example motivates the use of Frank-Wolfe algorithms in the context of generalized self-concordant functions.
We present more examples in the computational results.

\begin{example}[Intersection of a convex set with a polytope]\label{ex:intersection}
Consider Problem~\eqref{prob:main} where
$\mathcal{X} = \mathcal{P} \cap \mathcal{C}$, $\mathcal{P}$ is a polytope
over which we can minimize a linear function efficiently, and $\mathcal{C}$ is a convex compact set for which one can easily build a barrier function.

\begin{figure*}[ht!]
    \centering
    \hspace{\fill}
    \subfloat[Plot of $f(\vx)$.]{{\includegraphics[width= .3\textwidth]{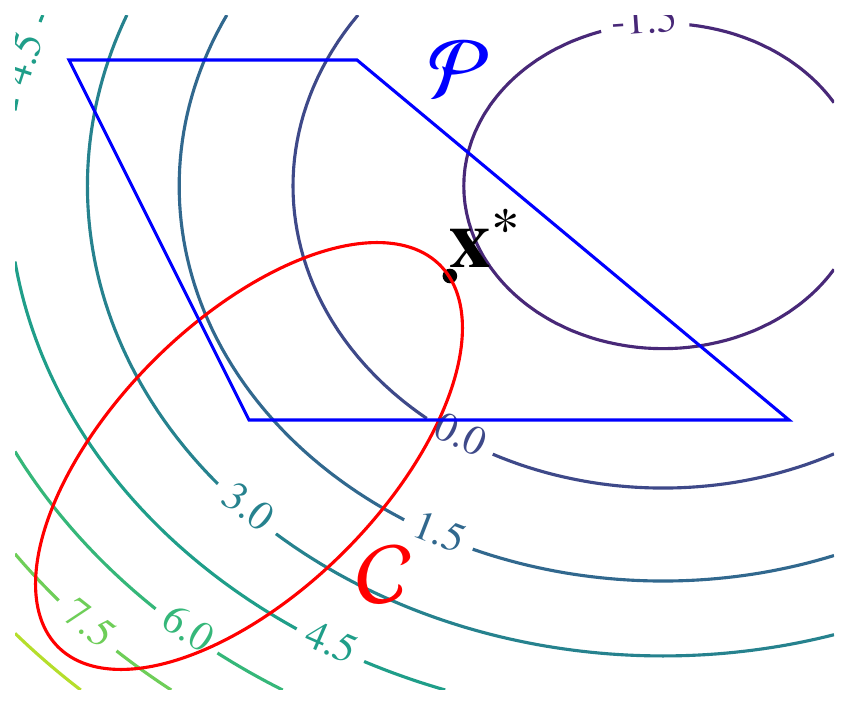} }\label{fig:no_log_barrier}}%
    \hspace{\fill}
    \subfloat[Plot of $f(\vx) + \mu' \Phi_{\mathcal{C}}(\vx)$.]{{\includegraphics[width = .3\textwidth]{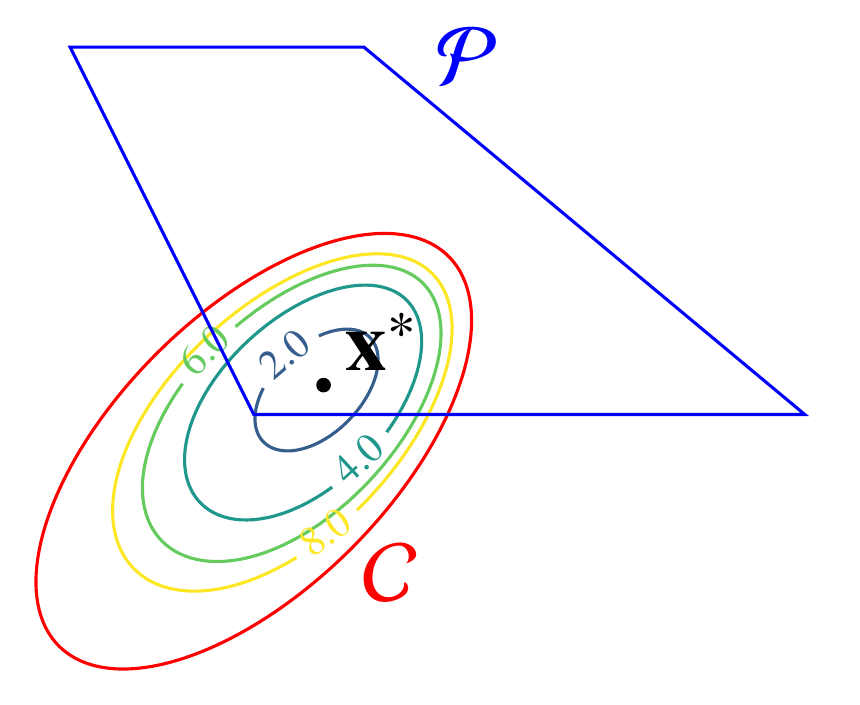} }\label{fig:log_barrier_low}}%
    \hspace{\fill}
    \subfloat[Plot of $f(\vx) + \mu \Phi_{\mathcal{C}}(\vx)$.]{{\includegraphics[width = .3\textwidth]{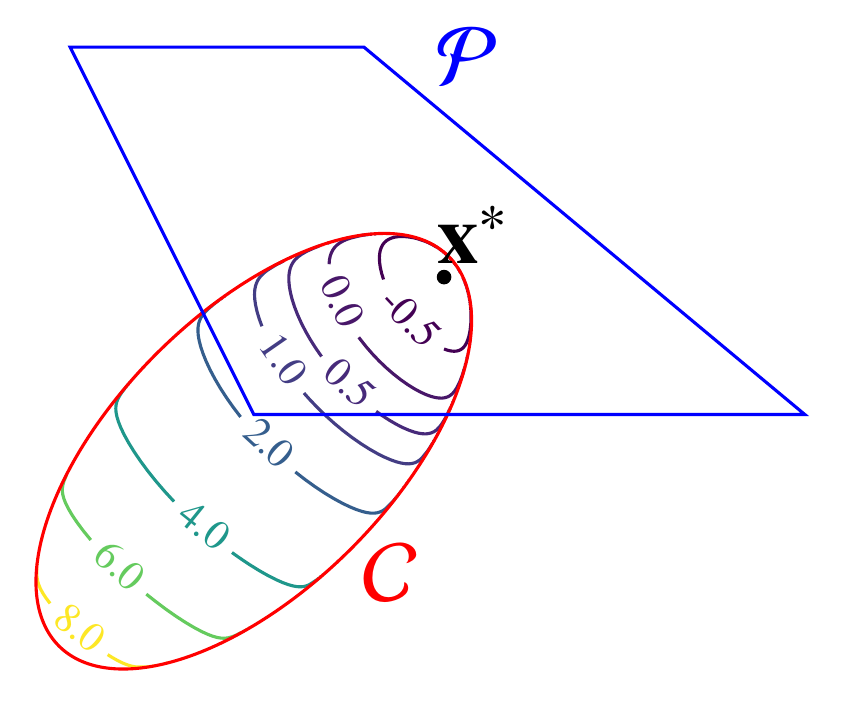} }\label{fig:log_barrier}}%
    \hspace*{\fill}
    \caption{Minimizing $f(\vx)$ over $\mathcal{P} \cap \mathcal{C}$, versus minimizing the sum of $f(\vx)$ and $\Phi_{\mathcal{C}}(\vx)$ over $\mathcal{P}$ for two different penalty values $\mu'$ and $\mu$ such that $\mu' \gg \mu$.}%
    \label{fig:problem_setting_schematic_convex}%
\end{figure*}

Solving a linear optimization problem over $\mathcal{X}$ may be extremely expensive. In light
of this, we can incorporate $\mathcal{C}$ into the problem through the use of a 
barrier penalty in the objective function, minimizing instead $f(\vx) + \mu \Phi_{\mathcal{C}}(\vx)$
where $\Phi_{\mathcal{C}}(\vx)$ is a log-barrier function for $\mathcal{C}$
and $\mu$ is a parameter controlling the penalization.
This reformulation is illustrated in Figure~\ref{fig:problem_setting_schematic_convex}.
Note that if the original objective function is generalized self-concordant, so is the new
objective function (see Proposition 1 in \cite{sun_generalized_2019}). We assume that computing the gradient of $f(\vx) + \mu \Phi_{\mathcal{C}}(\vx)$ is roughly as expensive as computing the gradient for $f(\vx)$ and solving an LP over $\mathcal{P}$ is inexpensive relative to solving an LP over $\mathcal{P}\cap \mathcal{C}$. The $\mu$ parameter can be driven down to $0$ after a solution converges
in a warm-starting procedure similar to interior-point methods,
ensuring convergence to the true optimum.

An additional advantage of this transformation of the problem is the solution structure.
Running Frank-Wolfe on the set $\mathcal{P} \cap \mathcal{C}$ can select a large number
of extremal points from $\boundary{\mathcal{C}}$ if $\mathcal{C}$ is non-polyhedral.
In contrast, $\mathcal{P}$ has a finite number of vertices, a small subset of which will be selected throughout the optimization procedure.
The same solution as that of the original problem can thus be constructed as a convex combination of
a small number of vertices of $\mathcal{P}$, improving sparsity and interpretability in many applications.
\end{example}

The following definition formalizes the setting of Problem~\eqref{prob:main}.

\begin{definition}[Generalized self-concordant function]
Let $f \in C^3\left( \dom f \right)$ be a closed convex function with $\dom f \subseteq \rr^n$ open. Then $f$ is $\left(M, \nu \right)$ \emph{generalized self-concordant} if:
\begin{align*}
|\innp{\operatorname{D}^3 f(\vx)[\vw]\vu, \vu}|\leq M \norm{\vu}_{\nabla^2 f(\vx)}^2 \norm{\vw}_{\nabla^2 f(\vx)}^{\nu - 2} \norm{\vw}_{2}^{3 - \nu},
\end{align*}
for any $\vx \in \dom f$ and $\vu,\vw \in \rr^n$, where
\begin{align*}
  \operatorname{D}^3 f(\vx)[\vw] = \lim\limits_{\alpha \rightarrow 0} \alpha^{-1}\left( \nabla^2 f(\vx + \alpha \vw) - \nabla^2 f(\vx)\right).
\end{align*}
\end{definition}
\vspace*{-.6cm}

\section{Frank-Wolfe Convergence Guarantees}

We establish convergence rates for a Frank-Wolfe variant with an open-loop step size strategy for generalized self-concordant functions. The \emph{Monotonic Frank-Wolfe (\texttt{M-FW})} algorithm presented in \cref{fw} is a rather simple, but powerful modification of the standard Frank-Wolfe algorithm, with the only difference that before taking a step, we verify if $\vx_t+\gamma_t(\mathbf{v}_t-\vx_t) \in \dom f$, and if so, we check whether moving to the next iterate provides primal progress.

\begin{algorithm}[H]
	\caption{Monotonic Frank-Wolfe (\texttt{M-FW})}
  \footnotesize
\label{fw}  
\begin{algorithmic}[1]
\Require Point $\vx_0\in \cx \cap \dom{f}$, function \(f\)
\hrulealg
\For{$t=0$ \textbf{to} $\dotsc$}
  \State$\mathbf{v}_t\leftarrow\argmin_{\mathbf{v}\in\cx}\innp{\nabla f(\vx_t),\mathbf{v}}$, $\gamma_t \leftarrow 2/(t+2)$
    \label{alg:fw_vertex}
  \State$\vx_{t+1}\leftarrow \vx_t+\gamma_t(\mathbf{v}_t-\vx_t)$ \label{alg:progress_is_made}
  \If{$\vx_{t+1} \notin \dom f$ \textbf{ or } $f(\vx_{t+1}) > f(\vx_t)$} \label{alg:domain:if}
  \State$\vx_{t+1}\leftarrow \vx_t$  \label{alg:same_point}
\EndIf
\EndFor
\end{algorithmic}
\end{algorithm}

Note, that the open-loop step size rule $2/(t+2)$ does not guarantee monotonic primal progress for the vanilla Frank-Wolfe algorithm in general.
If either of these two checks fails, we simply do not move: the algorithm sets $\vx_{t+1} = \vx_t$ in Line~\ref{alg:same_point} of \cref{fw}.
As customary, we assume short-circuit evaluation of the logical conditions in \cref{fw}, i.e., if the first condition in Line~\ref{alg:domain:if} is true, then the second condition is not even checked, and the algorithm directly goes to Line~\ref{alg:same_point}. This minor modification of the vanilla Frank-Wolfe algorithm enables us to use the monotonicity of the iterates in the proofs to come, at the expense of at most one extra function evaluation per iteration. Note that if we set $\vx_{t+1} = \vx_t$, we do not need to call the FOO or LMO oracle at iteration $t+1$, as we can simply reuse $\nabla f(\vx_t)$ and $\mathbf{v}_t$. This effectively means that between successive iterations in which we search for an acceptable value of $\gamma_t$, we only need to call the zeroth-order and domain oracle.

In order to establish the main convergence results for the algorithm, we lower bound the progress per iteration with the help of Proposition~\ref{proposition:key_inequalities}.

\begin{proposition}[Proposition 10, \cite{sun_generalized_2019}]\label{proposition:key_inequalities}
Given a $\left(M, \nu\right)$ generalized self-concordant function, then for $\nu\geq 2$, we have that:
\begin{align}
 f(\vy) - f(\vx) - \innp{\nabla f(\vx), \vy-\vx} \leq  \omega_{\nu}(d_{\nu}(\vx,\vy)) \norm{\vy - \vx}_{\nabla^2 f(\vx)}^2, \label{eq:pseudosmoothness_inequality}
\end{align}
where the inequality holds if and only if $d_{\nu}(\vx, \vy) < 1$ for $\nu > 2$, and we have that,
\begin{equation} \label{eq:distance_definition}
  d_{\nu}(\vx,\vy) \defeq \begin{cases}
    M \norm{\vy-\vx} & \mbox{if } \nu=2\\
    (\frac{\nu}{2} - 1) M \norm{\vy-\vx}^{3-\nu} \norm{\vy-\vx}_{\nabla^2 f(\vx)}^{\nu-2} & \mbox{if } \nu > 2,
  \end{cases}
\end{equation} 
\noindent
where:
\begin{equation*}
  \omega_\nu(\tau) \defeq \begin{cases}
    \frac{e^\tau - \tau - 1}{\tau^2} & \mbox{if } \nu = 2 \\
    \frac{-\tau - \text{ln}(1-\tau)}{\tau^2} & \mbox{if } \nu = 3 \\
    \frac{(1-\tau) \text{ln}(1-\tau) + \tau }{\tau^2} & \mbox{if } \nu = 4 \\
    \left(\frac{\nu-2}{4-\nu}\right) \,\,\frac{1}{\tau} \left[ \frac{\nu-2}{2(3-\nu)\tau} \left((1-\tau)^{\frac{2(3-\nu)}{2-\nu}} - 1 \right) -1 \right] & \mbox{otherwise.}
  \end{cases}
\end{equation*}
\end{proposition}

The inequality shown in \cref{eq:pseudosmoothness_inequality} is very similar to the one that we would obtain if the gradient of $f$ was Lipschitz continuous, however, while the Lipschitz continuity of the gradient leads to an inequality that holds globally for all $\vx,\vy \in \dom f$, the inequality in \eqref{eq:pseudosmoothness_inequality} only holds for $d_{\nu}(\vx,\vy) < 1$. Moreover, there are two other important differences, the norm used in \eqref{eq:pseudosmoothness_inequality} is now the norm defined by the Hessian at $\vx_t$ instead of the $\ell_2$ norm, and the term multiplying the norm is $\omega_{\nu}(d_{\nu}(\vx, \vy))$ instead of $1/2$. We deal with the latter issue by bounding $\omega_{\nu}(d_{\nu}(\vx, \vy))$ with a constant that depends on $\nu$ for any $\vx,\vy \in \dom f$ such that $d_{\nu}(\vx, \vy) \leq 1/2$, as shown in Remark~\ref{remark:bound_omega}.

\begin{remark}
\label{remark:bound_omega}
As $\operatorname{d}{\omega_{\nu}(\tau)}/\operatorname{d}{\tau} > 0$ for $\tau < 1$ and $\nu\geq 2$, then $\omega_{\nu}(\tau) \leq \omega_{\nu}(1/2)$ for $\tau \leq 1/2$.
\end{remark}

Due to the fact that we use a simple step size $\gamma_t = 2/(t+2)$, that we make monotonic progress, and we ensure that the iterates are inside $\dom f$, careful accounting allows us to bound the number of iterations until $d_{\nu}(\vx_{t}, \vx_t + \gamma_t(\mathbf{v}_t - \vx_t)) \leq 1/2$. Before formalizing the convergence rate, we first review a lemma needed in the proof.

\begin{lemma}[Proposition 7,\cite{sun_generalized_2019}]\label{lemma:inside_domain}
Let $f$ be a generalized self-concordant function with $\nu > 2$.
If $d_{\nu}(\vx,\vy) < 1$ and $\vx\in \dom f$ then $\vy \in \dom f$. For the case $\nu = 2$ we have that $\dom f = \rr^n$.
\end{lemma}

Putting all these things together allows us to obtain a convergence rate for \cref{fw}.

\begin{theorem}
\label{th:FW_generalized self-concordant_convergence}
Suppose $\cx$ is a compact convex set and $f$ is a $\left(M, \nu\right)$ generalized self-concordant function with $\nu \geq 2$, and define
the compact set
\begin{align*}
  \cl_0 \defeq \left\{ \vx \in \dom f \cap \cx \mid f(\vx) \leq f(\vx_0)\right\}.  
\end{align*}
Then, the Monotonic Frank-Wolfe algorithm (\cref{fw}) satisfies:
\begin{align}
h(\vx_t) \leq  \frac{4(T_{\nu}+1)}{t + 1}  \max \left\{h(\vx_0), L_f^{\cl_0}D^2 \omega_{\nu}(1/2) \right\}. \label{eq:convergence_rate}
\end{align}
for $t \geq T_{\nu}$, $T_{\nu}$ is defined as:
\begin{equation}
  T_{\nu} \defeq  \begin{cases}
    \left\lceil 4 MD \right\rceil - 2  & \mbox{if } \nu = 2 \\
     \left\lceil  2MD(L_f^{\cl_0})^{\nu/2-1}(\nu - 2)\right\rceil  - 2  & \mbox{otherwise}.
  \end{cases} \label{eq:burnin_period_value}
\end{equation}
Otherwise it holds that $h(\vx_t) \leq h(\vx_0)$ for $t<T_{\nu}$.
\begin{proof*}
As the algorithm makes monotonic progress and moves towards points such that $\vx_t \in \mathcal{X}$, then $\vx_t \in \cl_0$ for $t \geq 0$.
As the smoothness parameter of $f$ is bounded over $\cl_0$, we have from the properties of smooth functions that the bound $\norm{\vd}^2_{\nabla^2 f(\vx_t)} / \norm{\vd}^2 \leq L_f^{\cl_0}$ holds for any $\vd \in \rr^n$.
Particularizing for $\vd = \vx_t -\mathbf{v}_t$ and noting that $\norm{\vx_t -\mathbf{v}_t} \leq D$ leads to $\norm{\vx_t - \mathbf{v}_t}^2_{\nabla^2 f(\vx_t)} \leq L_f^{\cl_0} D^2$.
We then define $T_{\nu}$ with \eqref{eq:burnin_period_value}.
Using the definition shown in \eqref{eq:distance_definition} we have that for $t \geq T_{\nu}$ then 
$d(\vx_t, \vx_t + \gamma_t(\mathbf{v}_t - \vx_t)) \leq 1/2$. This fact, along with the fact that $\vx_t \in \dom f$ (by monotonicity) allows us to claim that $\vx_t + \gamma_t(\mathbf{v}_t - \vx_t) \in \dom f$, by application of Lemma~\ref{lemma:inside_domain}. This means that the non-zero step size $\gamma_t$ will ensure that $\vx_t + \gamma_t(\mathbf{v}_t - \vx_t) \in \dom f$ in Line~\ref{alg:domain:if} of \cref{fw}. Moreover, it allows us to use the bound between the function value at points $\vx_t$ and $\vx_{t} + \gamma_t \left( \mathbf{v}_t - \vx_t\right)$ in \eqref{eq:pseudosmoothness_inequality} of Proposition~\ref{proposition:key_inequalities}, which holds for $d(\vx_t, \vx_{t} + \gamma_t \left( \mathbf{v}_t - \vx_t\right)) < 1$. With this we can estimate the primal progress we can guarantee for $t \geq T_{\nu}$ if we move from $\vx_t$ to $\vx_{t} + \gamma_t \left( \mathbf{v}_t - \vx_t\right)$:
\begin{align*}
  h(\vx_{t+1}) &= h(\vx_{t} + \gamma_t \left( \mathbf{v}_t - \vx_t\right))\\
  &\leq h(\vx_t) - \gamma_t g(\vx_t) + \gamma_t^2 \omega_{\nu}\left( d_{\nu}\left( \vx_t, \vx_{t} + \gamma_t \left( \mathbf{v}_t - \vx_t\right) \right) \right)\norm{\mathbf{v}_t - \vx_t}^2_{\nabla^2 f(\vx_t)}\\
&\leq h(\vx_t)\left( 1 -  \gamma_t \right) + \gamma_t^2 L_f^{\cl_0} D^2 \omega_{\nu}(1/2), 
\end{align*}
where the second inequality follows from the upper bound on the primal gap via the FW gap $g(\vx_t)$, the application of Remark~\ref{remark:bound_omega} as for $t\geq T_{\nu}$ we have that $d_{\nu}(\vx_{t}, \vx_{t} + \gamma_t \left( \mathbf{v}_t - \vx_t\right)) \leq 1/2$, and from the fact that $\vx_t\in \cl_0$ for all $t\geq 0$. With the previous chain of inequalities we can bound the primal progress for $t \geq T_{\nu}$ as 
\begin{align}
h(\vx_t) - h(\vx_{t} + \gamma_t \left( \mathbf{v}_t - \vx_t\right)) & \geq  \gamma_t h(\vx_t) - \gamma_t^2L_f^{\cl_0}D^2\omega_{\nu}(1/2). \label{eq:primal_progress_key_ineq}
\end{align}
From these facts we can prove the convergence rate shown in \eqref{eq:convergence_rate} by induction. The base case $t = T_{\nu}$ holds trivially by the fact that using monotonicity we have that $h(\vx_{T_\nu}) \leq h(\vx_0)$. Assuming the claim is true for some $t \geq T_{\nu}$ we distinguish two cases. \\
\textbf{Case $\gamma_t h(\vx_t) - \gamma_t^2L_f^{\cl_0}D^2\omega_{\nu}(1/2)>0$: } Focusing on the first case, we can plug the previous inequality into \eqref{eq:primal_progress_key_ineq} to find that $\gamma_t$ guarantees primal progress, that is, $h(\vx_t) > h(\vx_{t} + \gamma_t \left( \mathbf{v}_t - \vx_t\right))$ with the step size $\gamma_t$, and so we know that we will not go into Line~\ref{alg:same_point} of \cref{fw}, and we have that $h(\vx_{t+1}) = h(\vx_{t} + \gamma_t \left( \mathbf{v}_t - \vx_t\right))$. Thus using the induction hypothesis and plugging in the expression for $\gamma_t = 2/(t+2)$ into \eqref{eq:primal_progress_key_ineq} we have:
\begin{align*}
h(\vx_{t+1}) & \leq   4\max \left\{h(\vx_0), L_f^{\cl_0}D^2 \omega_{\nu}(1/2)\right\}  \left( \frac{(T_{\nu}+1)t}{(t + 1)(t+2)} + \frac{1}{(t + 2)^2}  \right) \\
& \leq \frac{4(T_{\nu}+1)}{t + 2}\max \left\{h(\vx_0), L_f^{\cl_0}D^2 \omega_{\nu}(1/2) \right\},
\end{align*}
where we use that $(T_{\nu}+1)t/(t + 1) + 1/(t + 2)  \leq T_{\nu}+1$ for all $t\geq 0$ and any $t\geq T_{\nu}$. \\
\textbf{Case $\gamma_t h(\vx_t) - \gamma_t^2L_f^{\cl_0}D^2\omega_{\nu}(1/2) \leq 0$: } In this case, we cannot guarantee that the step size $\gamma_t$ provides primal progress by plugging into \eqref{eq:primal_progress_key_ineq}, and so we cannot guarantee if a step size of $\gamma_t$ will be accepted and we will have $\vx_{t+1} = \vx_{t} + \gamma_t \left( \mathbf{v}_t - \vx_t\right)$, or we will simply have $\vx_{t+1} = \vx_{t}$, that is, we may go into Line~\ref{alg:same_point} of \cref{fw}. Nevertheless, if we reorganize the expression $\gamma_t h(\vx_t) - \gamma_t^2L_f^{\cl_0}D^2\omega_{\nu}(1/2)\leq 0$, by monotonicity we will have that:
\begin{align*}
h(\vx_{t+1})  \leq h(\vx_t)  \leq \frac{2}{t + 2} L_f^{\cl_0}D^2\omega_{\nu}(1/2) \leq  \frac{4(T_{\nu} + 1)}{t + 2} \max \left\{h(\vx_0), L_f^{\cl_0}D^2 \omega_{\nu}(1/2) \right\}.
\end{align*}
Where the last inequality holds as $2 \leq 4(T_{\nu} + 1)$ for any $T_{\nu}\geq 0$.
\end{proof*}
\end{theorem}

One of the quantities that we have used in the proof of Theorem~\ref{th:FW_generalized self-concordant_convergence} is $L_f^{\cl_0}$. Note that the function $f$ is $L_f^{\cl_0}$-smooth over $\cl_0$. One could wonder why we have bothered to use the bound on the Bregman divergence in Proposition~\ref{proposition:key_inequalities} for a $(M, \nu)$-generalized self-concordant function, instead of simply using the bounds from the $L_f^{\cl_0}$-smoothness of $f$ over $\cl_0$. The reason is that the upper bound on the Bregman divergence in Proposition~\ref{proposition:key_inequalities} applies for any $\vx, \vy \in \dom f$ such that $d_{\nu}(\vx, \vy) < 1$, and we can easily bound the number of iterations $T_{\nu}$ it takes for the step size $\gamma_t = 2/(t+2)$ to verify both $\vx_t, \vx_{t} + \gamma_t (\mathbf{v}_t - \vx_t) \in \dom f$ and $d_{\nu}(\vx_t, \vx_{t} + \gamma_t (\mathbf{v}_t - \vx_t)) < 1$ for $t\geq T_{\nu}$. However, in order to apply the bound on the Bregman divergence from  $L_f^{\cl_0}$-smoothness we need $\vx_t, \vx_t + \gamma_t (\mathbf{v}_t - \vx_t) \in \cl_0$, and while it is easy to show by monotonicity that $\vx_t \in \cl_0$,  there is no straightforward way to prove that for some $\tilde{T}_{\nu}$ we have that $\vx_t + \gamma_t (\mathbf{v}_t - \vx_t)\in \cl_0$ for all $t \geq \tilde{T}_{\nu}$, i.e., that from some point onward a step with a non-zero step size is taken (that is, we do not go into Line~\ref{alg:same_point} of \cref{fw}) that guarantees primal progress.

\begin{remark}
In the case where $\nu = 2$ we can easily bound the primal gap $h(\vx_1)$, as in this setting $\dom f = \rr^n$, which leads to $h(\vx_1) \leq L_{f}^{\cx}D^2$ from \eqref{eq:primal_progress_key_ineq}, regardless of whether we set $\vx_{1} = \vx_0$ or $\vx_1 = \mathbf{v}_0$.
Moreover, as the upper bound on the Bregman divergence holds for $\nu = 2$ regardless of the value of $d_{2}(\vx, \vy)$, we can modify the proof of Theorem~\ref{th:FW_generalized self-concordant_convergence} to obtain a convergence rate of the form:
\begin{align*}
  h(\vx_{t}) \leq \frac{2}{t+1} L_f^{\cx}D^2 \omega_{2}(MD)\;\; \forall t \geq 1,
\end{align*}
which is reminiscient of the $\mathcal O(L_f^{\cx}D^2/t)$ rate of the original Frank-Wolfe algorithm for the smooth and convex case.
\end{remark}

Furthermore, with this simple step size we can also prove a convergence rate for the Frank-Wolfe gap, as shown in Theorem~\ref{th:FW_generalized self-concordant_convergence_FW_main_body}. More specifically, the minimum of the Frank-Wolfe gap over the run of the algorithm converges at a rate of $\mathcal{O}(1/t)$. The idea of the proof is very similar to the one in \cite{jaggi2013revisiting}. In a nutshell, as the primal progress per iteration is directly related to the step size times the Frank-Wolfe gap, we know that the Frank-Wolfe gap cannot remain indefinitely above a given value, as otherwise we would obtain a large amount of primal progress, which would make the primal gap become negative. This is formalized in Theorem~\ref{th:FW_generalized self-concordant_convergence_FW_main_body}.

\begin{theorem}
\label{th:FW_generalized self-concordant_convergence_FW_main_body}
Suppose $\cx$ is a compact convex set  and $f$ is a $\left(M, \nu\right)$ generalized self-concordant function with $\nu \geq 2$. Then if the Monotonic Frank-Wolfe algorithm (\cref{fw}) is run for $T \geq T_{\nu} + 6$ iterations, we will have that:
\begin{align*}
\min\limits_{1 \leq t \leq T} g(\vx_t) \leq \mathcal{O}(1/T),
\end{align*}
where $T_{\nu}$ is defined as: 
\begin{align}
  T_{\nu} \defeq  \begin{cases}
    \left\lceil 4 MD \right\rceil - 2  & \mbox{if } \nu = 2 \\
     \left\lceil  2MD(L_f^{\cl_0})^{\nu/2-1}(\nu - 2)\right\rceil  - 2  & \mbox{otherwise}.
  \end{cases} \label{eq:burnin_period_value_appx}
\end{align}
\begin{proof*}
In order to prove the claim, we focus on the iterations $t$ such that:
\begin{equation}\label{eq:tempineqonT}
  T_{\nu} + \lceil (T - T_{\nu})/3 \rceil - 2\leq t \leq T - 2,  
\end{equation}
where $T_{\nu}$ is defined in \eqref{eq:burnin_period_value_appx}.
Note that as we assume that $T \geq T_{\nu} + 6$, we know that $T_{\nu} \leq T_{\nu} + \lceil (T - T_{\nu})/3 \rceil - 2$, and so for iterations $T_{\nu} + \lceil (T - T_{\nu})/3 \rceil - 2\leq t \leq T - 2$ we know that $d_{\nu}(\vx_t, \vx_{t+1}) \leq 1/2$, and so:
\begin{align}
h(\vx_{t+1}) &\leq h(\vx_t) -  \gamma_t g(\vx_t)  + \gamma_t^2 L_f^{\cl_0} D^2 \omega_{\nu}(1/2). \label{Eq:contraction_general_nu_FW_gap}
\end{align}
In a very similar fashion as was done in the proof of Theorem~\ref{th:FW_generalized self-concordant_convergence}, we divide the proof into two different cases. \\
\textbf{Case $- \gamma_t g(\vx_t) + \gamma_t^2L_f^{\cl_0}D^2\omega_{\nu}(1/2) \geq 0$ for some $T_{\nu} + \lceil (T - T_{\nu})/3 \rceil - 2\leq t \leq T - 2$: }
Reordering the inequality above we therefore know that there exists a $T_{\nu} + \lceil (T - T_{\nu})/3 \rceil - 2\leq K \leq T - 2$ such that:
\begin{align*}
g(\vx_{K}) &\leq \frac{2}{2 + K} L_f^{\cl_0}D^2\omega_{\nu}(1/2) \\
&\leq \frac{2}{T_{\nu} + \lceil (T - T_{\nu})/3 \rceil} L_f^{\cl_0}D^2\omega_{\nu}(1/2) \\
& = \frac{6}{2T_{\nu} + T} L_f^{\cl_0}D^2\omega_{\nu}(1/2),
\end{align*}
where the second inequality follows from the fact that $T_{\nu} + \lceil (T - T_{\nu})/3 \rceil - 2\leq K$. This leads to $\min\limits_{1 \leq t \leq T} g(\vx_t) \leq g(\vx_{K}) \leq \frac{6}{ 2T_{\nu} + T} L_f^{\cl_0}D^2\omega_{\nu}(1/2)$. \\
\textbf{Case $- \gamma_t g(\vx_t) + \gamma_t^2L_f^{\cl_0}D^2\omega_{\nu}(1/2) < 0$ for all $T_{\nu} + \lceil (T - T_{\nu})/3 \rceil - 2\leq t \leq T - 2$: } 
Using the inequality above and plugging into \eqref{Eq:contraction_general_nu_FW_gap} allows us to conclude that all steps $T_{\nu} + \lceil (T - T_{\nu})/3 \rceil- 2\leq t \leq T - 2$ will produce primal progress using the step size $\gamma_t$, and so as we know that $\vx_{t + 1}\in \dom f$ by Lemma~\ref{lemma:inside_domain}, then for all $T_{\nu} + \lceil (T - T_{\nu})/3 \rceil - 2\leq t \leq T - 2$ we will take a non-zero step size determined by $\gamma_t$, as $\vx_{t} + \gamma_t (\mathbf{v}_t - \vx_t)\in \dom f$ and $f(\vx_{t} + \gamma_t (\mathbf{v}_t - \vx_t)) < f(\vx_t)$ in \cref{alg:domain:if} of \cref{fw}. Consequently, summing up \eqref{Eq:contraction_general_nu_FW_gap} from $t_{\min} \defeq T_{\nu} + \lceil (T - T_{\nu})/3 \rceil - 2$ to $t_{\max} \defeq T - 2$ we have that:
\begin{align}
h(\vx_{t_{\max} + 1}) &\leq  h\left(\vx_{t_{\min}}\right) -  \sum_{t = t_{\min}}^{t_{\max}} \gamma_t g(\vx_t) + L_f^{\cl_0} D^2 \omega_{\nu}(1/2) \sum_{t = t_{\min}}^{t_{\max}}  \gamma_t^2 \\
& \leq   h\left(\vx_{t_{\min}}\right) -  2\min\limits_{ t_{\min} \leq t \leq t_{\max}} g(\vx_t) \sum_{t = t_{\min}}^{t_{\max}} \frac{1}{2 + t} + 4L_f^{\cl_0} D^2 \omega_{\nu}(1/2) \sum_{t = t_{\min}}^{t_{\max}} \frac{1}{(2 + t)^2} \label{eq:take_minimum} \\
& \leq   h\left(\vx_{t_{\min}}\right) -  2\min\limits_{ 1 \leq t \leq T} g(\vx_t) \frac{t_{\max} - t_{\min} + 1}{2 + t_{\max}} + 4L_f^{\cl_0} D^2 \omega_{\nu}(1/2)  \frac{t_{\max} - t_{\min} + 1}{(2 + t_{\min})^2} \label{eq:bounds_iteration} \\
& \leq   4\left( \frac{T_{\nu}+1}{t_{\min} + 1}  + \frac{t_{\max} - t_{\min} + 1}{(2 + t_{\min})^2} \right) \max \left\{h(\vx_0), L_f^{\cl_0}D^2 \omega_{\nu}(1/2) \right\} \label{eq:plug_primal_gap1}\\
&\;\;\;\;  -  2\min\limits_{ 1 \leq t \leq T} g(\vx_t) \frac{t_{\max} - t_{\min} + 1}{2 + t_{\max}}. \label{eq:plug_primal_gap2}
\end{align}
Note that \eqref{eq:take_minimum} stems from the fact that $\min_{ t_{\min} \leq t \leq t_{\max}} g(\vx_t) \leq g(\vx_t)$ for any $t_{\min} \leq t \leq  t_{\max}$, and from plugging $\gamma_t = 2/(2+t)$, and \eqref{eq:bounds_iteration} follows from the fact that $-1/(2+t)\leq -1/(2 + t_{\max})$ and $1/(2+t)\leq 1/(2 + t_{\min})$ for all $t_{\min} \leq t \leq  t_{\max}$. The last inequality, shown in \eqref{eq:plug_primal_gap1} and \eqref{eq:plug_primal_gap2} arises from plugging in the upper bound on the primal gap $h(\vx_{t_{\min}})$ from Theorem~\ref{th:FW_generalized self-concordant_convergence} and collecting terms. If we plug in the specific values of $t_{\max}$ and $t_{\min}$ this leads to:
\begin{align}
h(\vx_{T - 1}) & \leq   12 \left( \frac{T_{\nu}+1}{2T_{\nu} + T- 3}  + \frac{2T - 2T_{\nu} + 3}{(2T_{\nu} + T)^2} \right) \max \left\{h(\vx_0), L_f^{\cl_0}D^2 \omega_{\nu}(1/2) \right\} \\
& -  \frac{2}{3}\min\limits_{ 1 \leq t \leq T} g(\vx_t) \frac{T - T_{\nu}}{T}. \label{eq:key_inequality}
\end{align}
We establish our claim using proof by contradiction. Assuming that:
\begin{align*}
\min\limits_{ 1 \leq t \leq T} g(\vx_t) > \frac{18 T}{T - T_{\nu}} \left( \frac{T_{\nu}+1}{2T_{\nu} + T- 3}  + \frac{2T - 2T_{\nu} + 3}{(2T_{\nu} + T)^2} \right) \max \left\{h(\vx_0), L_f^{\cl_0}D^2 \omega_{\nu}(1/2) \right\}
\end{align*}
results, together with the bound from \eqref{eq:key_inequality} in $h(\vx_{T - 1})<0$, which is the desired contradiction, as the primal gap cannot be negative.
Therefore we must have that:
\begin{align*}
\min\limits_{ 1 \leq i \leq T} g(\vx_i) & \leq \frac{18 T}{T - T_{\nu}} \left( \frac{T_{\nu}+1}{2T_{\nu} + T- 3}  + \frac{2T - 2T_{\nu} + 3}{(2T_{\nu} + T)^2} \right) \max \left\{h(\vx_0), L_f^{\cl_0}D^2 \omega_{\nu}(1/2) \right\}\\
& = \mathcal{O}(1/T).
\end{align*}
This completes the proof.
\end{proof*}
\end{theorem}

\begin{remark} 
Note that the Monotonic Frank-Wolfe algorithm (\cref{fw}) performs at most one ZOO, FOO, DO, and LMO oracle call per iteration.
This means that Theorems~\ref{th:FW_generalized self-concordant_convergence} and \ref{th:FW_generalized self-concordant_convergence_FW_main_body} effectively bound the number of ZOO, FOO, DO, and LMO oracle calls needed to achieve a target primal gap or Frank-Wolfe gap accuracy $\varepsilon$ as a function of $T_{\nu}$ and $\varepsilon$; note that $T_{\nu}$ is independent of $\varepsilon$. This is an important difference with respect to existing bounds, as the existing Frank-Wolfe-style first-order algorithms for generalized self-concordant functions in the literature that utilize various types of line searches may perform more than one ZOO or DO call per iteration in the line search. This means that the convergence bounds in terms of iteration count of these algorithms are only informative when considering the number of FOO and LMO calls that are needed to reach a target accuracy in primal gap, and do not directly provide any information regarding the number of ZOO or DO calls that are needed. In order to bound the latter two quantities one typically needs additional technical tools. For example, for the backtracking line search of  \cite{pedregosa2018step}, one can use \cite[Theorem 1, Appendix C]{pedregosa2018step}, or a slightly modified version of \cite[Lemma 4]{nesterov2013gradient}, to find a bound for the number of ZOO or DO calls that are needed to find an $\varepsilon$-optimal solution. Note that these bounds depend on user-defined initialization or tuning parameters.
\end{remark}

\begin{remark} \label{remark:halving}
In practice, a halving strategy for the step size is preferred for the 
implementation of the Monotonic Frank-Wolfe algorithm, as opposed to the step size implementation shown in \cref{fw}. This halving strategy, which is shown in \cref{alg:halvingmfw}, helps
deal with the case in which a large number of \textbf{consecutive} step sizes $\gamma_t$ 
are rejected either because $\vx_{t} + \gamma_t(\mathbf{v}_t - \vx_t) \notin \dom f$ or $f(\vx_{t} + \gamma_t(\mathbf{v}_t - \vx_t)) > f(\vx_t)$,
and helps avoid the need to potentially call the zeroth-order or domain oracle a large number
of times in these cases. The halving strategy in \cref{alg:halvingmfw} results in a step size that is at most a factor 
of $2$ smaller than the one that would have been accepted with the original strategy, i.e., that would have 
ensured that $\vx_{t} + \gamma_t( \mathbf{v}_t -\vx_{t}) \in \dom f$ and 
$f(\vx_{t} + \gamma_t( \mathbf{v}_t -\vx_{t})) \leq f(\vx_t)$, in the standard
Monotonic Frank-Wolfe algorithm in \cref{fw}. However, the number
of zeroth-order or domain oracles that would be needed to find this step size
that satisfies both $\vx_{t} + \gamma_t( \mathbf{v}_t -\vx_{t}) \in \dom f$ and
$f(\vx_{t+1}) \leq f(\vx_t)$ is logarithmic for the Monotonic Frank-Wolfe variant 
shown in \cref{alg:halvingmfw}, when compared to the number needed for
the Monotonic Frank-Wolfe variant without halving shown in \cref{fw}.
Note that the convergence properties established throughout the paper for the Monotonic Frank-Wolfe
algorithm in \cref{fw} also hold for the variant in \cref{alg:halvingmfw}; with the only difference being that we lose a very small constant factor (e.g., at most a factor of $2$ for the standard case) in
the convergence rate.
\end{remark}

\begin{algorithm}[H]
	\caption{Halving \texttt{M-FW}}
  \label{alg:halvingmfw}
  \footnotesize
\begin{algorithmic}[1]
\Require Point $\vx_0\in \cx \cap \dom{f}$, function \(f\)
\Ensure Iterates $\vx_1, \dotsc \in \cx$
\hrulealg
\State $\psi_{-1} \gets 0$
\For{$t=0$ \textbf{to} $\dotsc$}
  \State$\mathbf{v}_t\leftarrow\argmin_{\mathbf{v}\in\cx}\innp{\nabla f(\vx_t),\mathbf{v}}$
  \State $\psi_t \gets \psi_{t-1}$
  \State$\gamma_t \leftarrow 2^{1-\psi_t} /(t+2)$
  \State$\vx_{t+1}\leftarrow \vx_t+\gamma_t(\mathbf{v}_t-\vx_t)$
  \While{$\vx_{t+1} \notin \dom f$ \textbf{ or } $f(\vx_{t+1}) > f(\vx_t)$}
  \State $\psi_t \gets \psi_t + 1$
  \State$\gamma_t \leftarrow 2^{1-\psi_t} /(t+2)$
  \State$\vx_{t+1}\leftarrow \vx_t+\gamma_t(\mathbf{v}_t-\vx_t)$
\EndWhile
\EndFor
\end{algorithmic}
\end{algorithm}

In Table~\ref{table:bound_comparison} we provide a detailed complexity comparison between the Monotonic Frank-Wolfe (\texttt{M-FW}) algorithm (\cref{fw}), and other comparable algorithms in the literature.

\begin{table}[H]\centering
  \footnotesize
\begin{tabular}{lccccc}
\toprule
\textbf{Algorithm} & \textbf{SOO calls} & \textbf{FOO calls} & \textbf{ZOO calls}& \textbf{LMO calls} & \textbf{DO calls} \\
\midrule
\texttt{FW-GSC} \cite[Alg.2]{dvurechensky2020generalized}  & $\mathcal{O}(1/\varepsilon)$ & $\mathcal{O}(1/\varepsilon)$ &  & $\mathcal{O}(1/\varepsilon)$  &   \\
\texttt{LBTFW-GSC}$^{\ddagger}$ \cite[Alg.3]{dvurechensky2020generalized} &  & $\mathcal{O}(1/\varepsilon)$ & $\mathcal{O}(1/\varepsilon)$ & $\mathcal{O}(1/\varepsilon)$  & $\mathcal{O}(1/\varepsilon)$  \\
\texttt{MBTFW-GSC}$^{\ddagger}$ \cite[Alg.5]{dvurechensky2020generalized} & $\mathcal{O}(1/\varepsilon)$ & $\mathcal{O}(1/\varepsilon)$ & $\mathcal{O}(1/\varepsilon)$ & $\mathcal{O}(1/\varepsilon)$  & $\mathcal{O}(1/\varepsilon)$  \\
\texttt{M-FW}$^{\dagger}$ [\textbf{This work}] &  & $\mathcal{O}(1/\varepsilon)$ & $\mathcal{O}(1/\varepsilon)$ & $\mathcal{O}(1/\varepsilon)$ & $\mathcal{O}(1/\varepsilon)$  \\
\bottomrule
\end{tabular}
\smallskip
\caption{
\textbf{Complexity comparison}: Number of iterations needed to reach a solution with $h(\vx)$ below $\varepsilon$ for Problem~\ref{prob:main}. We use the superscript $\dagger$ to indicate that the same complexities hold for reaching an $\varepsilon$-optimal solution in $g(\vx)$. The superscript $\ddagger$ is used to indicate that constants in the convergence bounds depend on user-defined inputs; the other algorithms are parameter-free.}
 \label{table:bound_comparison}
\end{table}

We note that the \texttt{LBTFW-GSC} algorithm from \citet{dvurechensky2020generalized} is in essence the Frank-Wolfe algorithm with a modified version of the backtracking line search of \citet{pedregosa2018step}. In the next section, we provide improved convergence guarantees for various cases of interest for this algorithm, which we refer to as the Frank-Wolfe algorithm with Backtrack (\texttt{B-FW}) for simplicity.

\subsection{Improved convergence guarantees}

We will now establish improved convergence rates for various special cases. We focus on two different settings to obtain improved convergence rates;
in the first, we assume that $\vx^* \in \interior \left( \cx \cap \dom f \right)$ (Section~\ref{Appx:Optininterior}),
and in the second we assume that $\cx$ is strongly or uniformly convex (Section~\ref{Appx:str_cvx_sets}). The algorithm in this section is a slightly modified Frank-Wolfe algorithm with the adaptive line search technique of \citet{pedregosa2018step} (shown for reference in \cref{fw2} and \ref{backtrack}).
This is the same algorithm used in \citet{dvurechensky2020generalized}, however, we show improved convergence rates in several settings of interest.
Note that the adaptive line search technique of \citet{pedregosa2018step} requires user-defined inputs or parameters, which means that the algorithms in this section are not parameter-free.
The parameter $M$ of \cref{backtrack} corresponds to a local estimate of the Lipschitz constant of $f$, the stopping condition defining the admissible step size requires the function decrease to be greater than the one derived from the quadratic model built from the Lipschitz estimate $M$ and gradient, hence ensuring monotonicity.

\begin{algorithm}[H]
	\caption{\texttt{B-FW}}
  \label{fw2}
  \footnotesize
\begin{algorithmic}[1]
\Require Point $\vx_0\in \cx \cap \dom{f}$, function \(f\), initial smoothness estimate $L_{-1}$
\Ensure Iterates $\vx_1, \dotsc \in \cx$
\hrulealg
\For{$t=0$ \textbf{to} $\dotsc$}
  \State$\mathbf{v}_t\leftarrow\argmin_{\mathbf{v}\in\cx}\innp{\nabla f(\vx_t),\mathbf{v}}$
    \label{alg:fw_vertex2}
  \State $\gamma_t, L_{t} \leftarrow$ \texttt{Backtrack}$\left(f, \vx_t, \mathbf{v}_t - \vx_t, L_{t-1}, 1 \right)$
  \State$\vx_{t+1}\leftarrow \vx_t+\gamma_t(\mathbf{v}_t-\vx_t)$ \label{alg:progress_is_made2}
\EndFor
\end{algorithmic}
\end{algorithm}

\begin{algorithm}[H]
	\caption{\texttt{Backtrack}$\left(f, \vx, \vd, L_{t-1}, \gamma_{\max} \right)$ (line search of \cite{pedregosa2018step})}
\label{backtrack}
\footnotesize
\begin{algorithmic}[1]
\Require Point $\vx \in \cx \cap \dom{f}$, $\mathbf{v} \in \rr^n$, function $f$, estimate $L_{t-1}$, step size $\gamma_{\max}$
\Ensure $\gamma$, $M$
\hrulealg
\State Choose $\tau>1$, $\eta\leq 1$ and $M\in [\eta L_{t-1}, L_{t-1}]$
\State $\gamma = \min\{ -\innp{\nabla f(\vx), \vd}/(M\norm{\vd}^2), \gamma_{\max}\}$
\While{$\vx + \gamma\vd\notin \dom f$ \textbf{ or } $f(\vx + \gamma\vd) - f(\vx)> \frac{M\gamma^2}{2}  \norm{\vd}^2 + \gamma \innp{\nabla f(\vx), \vd}$} \label{backtrack:domTestBacktrack}
\State $M = \tau M$
\State $\gamma = \min\{ -\innp{\nabla f(\vx), \vd}/(M\norm{\vd}^2), \gamma_{\max}\}$ \label{backtrack:stepsize}
\EndWhile
\end{algorithmic}
\end{algorithm}

\subsubsection{Optimum contained in the interior}\label{Appx:Optininterior}

We first focus on the assumption that $\vx^* \in \interior \left(\cx \cap \dom f \right) $, obtaining improved rates when we use the FW algorithm coupled with the adaptive step size strategy from \citet{pedregosa2018step} (see \cref{backtrack}).
This assumption is reasonable if for example $\boundary{\cx} \not\subseteq \dom f$, and $\interior\left(\cx\right) \subseteq \dom f$.
That is to say, we will have that $\vx^* \in \interior \left(\cx \cap \dom f \right)$ if for example we use logarithmic barrier functions to encode a set of constraints, and we have that $\dom f$ is a proper subset of $\cx$.
In this case the optimum is guaranteed to be in $\interior \left(\cx \cap \dom f \right)$.

The analysis in this case is reminiscent of the one in the seminal work of \citet{guelat1986some}, and is presented in Subsection~\ref{Appx:Optininterior}.
Note that we can upper-bound the value of $L_t$ for $t\geq 0$ by $\tilde{L} \defeq \max\{ \tau L_f^{\cl_0}, L_{-1}\}$, where $\tau>1$ is the backtracking parameter and 
$L_{-1}$ is the initial smoothness estimate in \cref{backtrack}.
Before proving the main theoretical results of this section, we first review some auxiliary results that allow us to prove linear convergence in this setting.

\begin{proposition}[Proposition 3, \cite{sun_generalized_2019}] \label{prop:PD-Hessian}
Let $f$ be generalized self-concordant with $\nu \geq 2$ and $\dom f$ not contain any straight line, then the Hessian $\nabla^2 f(\vx)$ is non-degenerate at all points $\vx \in \dom f$.
\end{proposition}

Note that the assumption that $\dom{f}$ does not contain any straight line is
without loss of generality as we can simply modify the function outside of our
compact convex feasible region so that it holds.

\begin{proposition}[Proposition 2.16, \cite{cg_survey}] \label{th:bound_interior_optimum}
If there exists an $r> 0$ such that $\mathcal{B}(\vx^*, r)\subseteq \cx \cap \dom f$, then for all $\vx\in \cx \cap \dom f$ we have that:
\begin{align*}
\frac{g(\vx)}{\norm{\vx - \mathbf{v}}} \geq \frac{r}{D} \norm{\nabla f(\vx)} \geq \frac{r}{D}\frac{\innp{\nabla f(\vx), \vx - \vx^*}}{\norm{\vx - \vx^*}},
\end{align*}
where $\mathbf{v} = \argmin\limits_{\vy \in \cx} \innp{\nabla f(\vx), \vy}$ and $g(\vx)$ is the Frank-Wolfe gap.
\end{proposition}

With these tools at hand, we show that the Frank-Wolfe algorithm with the backtracking step-size strategy converges at a linear rate.

\begin{theorem}
\label{th:FW_generalizedself-concordant_convergence_str_cvx_case_inner_opt}
Let $f$ be a $\left(M, \nu\right)$ generalized self-concordant function with $\nu \geq 2$, let $\dom f$ not contain any straight line
and define the compact set
\begin{align*}
\cl_0 \defeq \left\{ \vx \in \dom f \cap \cx \mid f(\vx) \leq f(\vx_0)\right\}.  
\end{align*}
Furthermore, let $r > 0$ be the largest value such that $\mathcal{B}(\vx^*, r)\subseteq \cx \cap \dom f$.
Then, the Frank-Wolfe algorithm (\cref{fw2}) with the backtracking strategy of \citet{pedregosa2018step} results in a linear primal gap convergence rate of the form:
\begin{align*}
h(\vx_{t}) & \leq h(\vx_0) \left( 1 -  \frac{\mu_f^{\cl_0}}{2\tilde{L}}  \left(\frac{r}{D}\right)^2   \right)^t ,
\end{align*}
for $t\geq 1$, where $\tilde{L} \defeq \max\{ \tau L_f^{\cl_0}, L_{-1}\}$, $\tau>1$ is the backtracking parameter, 
$L_{-1}$ is the initial smoothness estimate in \cref{backtrack}.
\begin{proof}
As the backtracking line search makes monotonic primal progress, we know that for $t\geq 0$ we will have that $\vx_t \in \cl_0$.
Since $\dom f$ does not contain any straight line by assumption, we know from \cref{prop:PD-Hessian} that for all $\vx\in \dom f$, the Hessian is non-degenerate and therefore $\mu_f^{\cl_0}>0$. This allows us to claim that for any $\vx, \vy\in \cl_0$ we have that:
\begin{align}
f(\vx) - f(\vy) - \innp{\nabla f(\vy), \vx -\vy} &\geq \frac{\mu_f^{\cl_0}}{2}  \norm{\vx -\vy}^2. \label{appx:lower_bound_Bregman1}
\end{align}
The backtracking line search in \cref{backtrack} will either output a point $\gamma_t = 1$ or $\gamma_t < 1$.
In any case, \cref{backtrack} will find and output a smoothness estimate $L_t$ and a step size $\gamma_t$ such that for $\vx_{t+1} = \vx_t + \gamma_t(\mathbf{v}_t - \vx_t)$ we have that: 
\begin{align}
f(\vx_{t+1}) - f(\vx_t) \leq \frac{L_t\gamma_t^2}{2}  \norm{\vx_t - \mathbf{v}_t}^2 - \gamma_t g(\vx_t). \label{local:smoothness}
\end{align}
In the case where $\gamma_t = 1$ we know by observing Line~\ref{backtrack:stepsize} of \cref{backtrack} that $g(\vx_t) \geq L_t \norm{\vx_t - \mathbf{v}_t}^2$, and so plugging into \eqref{local:smoothness} we arrive at $h(\vx_{t+1}) \leq h(\vx_t)/2$.
In the case where $\gamma_t = g(\vx_t)/(L_t \norm{\vx_t - \mathbf{v}_t}^2 ) < 1$, we have that $g(\vx_t) < L_t \norm{\vx_t - \mathbf{v}_t}^2$, which leads to $h(\vx_{t+1}) \leq h(\vx_k) - g(\vx_t)^2/(2L_t \norm{\vx_t - \mathbf{v}_t}^2)$, when plugging the expression for the step size in the progress bound in \eqref{local:smoothness}.
In this last case where $\gamma_t < 1$ we have the following contraction for the primal gap:
\begin{align*}
h(\vx_t) - h(\vx_{t+1}) & \geq  \frac{g(\vx_t)^2}{2L_t\norm{\vx_t - \mathbf{v}_t}^2} \\
& \geq  \frac{r^2}{D^2}\frac{\norm{\nabla f(\vx_t)}^2}{2L_t} \\
& \geq  \frac{\mu_f^{\cl_0}}{\tilde{L}} \frac{r^2}{D^2}h(\vx_t),
\end{align*}
where we have used the inequality that involves the central term and the leftmost term in Proposition~\ref{th:bound_interior_optimum}, and the last inequality stems from the bound $h(\vx_t) \leq \norm{\nabla f(\vx_t)}^2/(2\mu_f^{\cl_0})$ for $\mu_f^{\cl_0}$-strongly convex functions. Putting the above bounds together we have that:
\begin{align*}
h(\vx_{t+1}) & \leq h(\vx_t) \left( 1 - \frac{1}{2} \min\left\{ 1,\frac{\mu_f^{\cl_0}}{\tilde{L}}  \left(\frac{r}{D}\right)^2  \right\} \right) \\
& \leq h(\vx_t) \left( 1 -  \frac{\mu_f^{\cl_0}}{2\tilde{L}}  \left(\frac{r}{D}\right)^2  \right),
\end{align*}
which completes the proof.
\end{proof}
\end{theorem}

The previous bound depends on the largest positive $r$ such that $\mathcal{B}(\vx^*, r)\subseteq \cx \cap \dom f$, which can be arbitrarily small.
Note also that the previous proof uses the lower bound of the Bregman divergence from the $\mu_{f}^{\cl_0}$-strong convexity of the function over $\cl_0$ to obtain linear convergence.
Note that this bound is local as this $\mu_{f}^{\cl_0}$-strong convexity holds only inside $\cl_0$, and is only of use because the step size of \cref{backtrack} automatically ensures that if $\vx_t \in \cl_0$ and $\vd_t$ is a descent direction, then $\vx_t + \gamma_t \vd_t \in \cl_0$.
This is in contrast with \cref{fw}, in which the step size $\gamma_t = 2/(2+t)$ did not automatically ensure monotonicity in primal gap, and this had to be enforced by setting $\vx_{t+1} = \vx_t$ if $f(\vx_t + \gamma_t \vd_t) > f(\vx_t)$, where $\vd_t = \mathbf{v}_t - \vx_t$.
If we were to have used the lower bound on the Bregman divergence from \cite[Proposition 10]{sun_generalized_2019} in the proof, which states that:
\begin{align*}
f(\vy) - f(\vx) - \innp{\nabla f(\vx), \vy-\vx} \geq \omega_{\nu}(-d_{\nu}(\vx,\vy)) \norm{\vy - \vx}_{\nabla^2 f(\vx)}^2, 
\end{align*}
for any $\vx, \vy \in \dom f$ and any $\nu \geq 2$, we would have arrived at a bound that holds over all $\dom f$. However, in order to arrive at a usable bound, and armed only with the knowledge that the Hessian is non-degenerate if $\dom f$ does not contain any straight line, and that $\vx, \vy \in \cl_0$, we would have had to write:
\begin{align*}
 \omega_{\nu}(-d_{\nu}(\vx,\vy)) \norm{\vx- \vy}_{\nabla^2 f(\vy)}^2 \geq  \mu_f^{\cl_0} \omega_{\nu}(-d_{\nu}(\vx,\vy))  \norm{\vx-\vy}^2,
\end{align*}
where the inequality follows from the definition of $\mu_{f}^{\cl_0}$. It is easy to see that as $d \omega_{\nu}(\tau)/d \tau > 0$ by Remark~\ref{remark:bound_omega}, we have that $1/2 =\omega_{\nu}(0) \geq \omega_{\nu}(-d_{\nu}(\vx,\vy))$. This results in a bound:
\begin{align}
f(\vy) - f(\vx) - \innp{\nabla f(\vx), \vy-\vx} \geq \mu_f^{\cl_0} \omega_{\nu}(-d_{\nu}(\vx,\vy))  \norm{\vx-\vy}^2. \label{appx:lower_bound_Bregman2}
\end{align}
When we compare the bounds obtained from local strong convexity in \eqref{appx:lower_bound_Bregman1} and that obtained directly from generalized self-concordance in \eqref{appx:lower_bound_Bregman2}, we can see that the former is tighter than the latter, albeit local.
For this reason, we have used the former bound in the proof of Theorem~\ref{th:FW_generalizedself-concordant_convergence_str_cvx_case_inner_opt}.

\subsubsection{Strongly convex or uniformly convex sets} \label{Appx:str_cvx_sets}

Next, we recall the definition of uniformly convex sets, used in \citet{kerdreux2021projection}, which will allow us to obtain improved convergence rates for the FW algorithm over uniformly convex feasible regions.

In order to prove convergence rate results for the case where the feasible region is $(\kappa, p)$-uniformly convex, we first review the definition of the $(\kappa, p)$-uniform convexity of a set (see \cref{appx:defn:uniform-convex-set}), as well as a useful lemma that allows us to go from contractions to convergence rates.

\begin{definition}[$(\kappa, q)$-uniformly convex set]
  \label{appx:defn:uniform-convex-set}
  Given two positive numbers $\kappa$ and $q$, we say the set $\cx\subseteq \rr^n$ is $(\kappa, q)$-\emph{uniformly convex} with respect to a norm $\norm{\cdot}$ if for any $\vx, \vy \in \cx$, $0\leq \gamma \leq 1$, and $\vz \in \rr^n$ with $\norm{\vz} = 1$ we have that: 
  \begin{equation*}
   \vy + \gamma (\vx - \vy) + \gamma(1 - \gamma) \cdot \kappa \norm{\vx - \vy}^q \vz \in \cx.
  \end{equation*}
\end{definition}

The previous definition allows us to obtain a scaling inequality very similar to the one shown in \cref{th:bound_interior_optimum}, which is key to proving the following convergence rates, and can be implicitly found in \cite{kerdreux2021projection} and \cite{garber2016linearly}.

\begin{proposition} \label{appx:key_proposition_str_cvx_sets}
Let $\cx\subseteq \rr^n$ be $(\kappa, q)$-\emph{uniformly convex}, then for all $\vx \in \cx$:
\begin{align*}
\frac{g(\vx)}{\norm{\vx - \mathbf{v}}^{q}} \geq \kappa \norm{\nabla f(\vx)},
\end{align*}
where $\mathbf{v} = \argmin_{\vu\in \cx}\innp{\nabla f(\vx), \vu}$, and $g(\vx)$ is the Frank-Wolfe gap.
\end{proposition}

The next lemma that will be presented is an extension of the one used in \cite[Lemma A.1]{kerdreux2021projection} (see also \citet{temlyakov2015greedy}), and allows us to go from per-iteration contractions to convergence rates.

\begin{lemma} \label{appx:contraction_to_convergence_rate}
We denote a sequence of nonnegative numbers by $\{ h_t\}_t$. Let $c_0$, $c_1$, $c_2$ and $\alpha$ be positive numbers such that $c_1 < 1$, $h_1 \leq c_0$ and $h_t - h_{t+1} \geq h_t \min\{c_1, c_2 h_t^{\alpha}\}$ for $t \geq 1$, then:
\begin{equation*}
  h_t \leq  \begin{cases}
     c_0 \left( 1 - c_1 \right)^{t-1}   & \mbox{if } 1 \leq t \leq t_0 \\
     \frac{\left( c_1/c_2\right)^{1/\alpha}}{\left(1 + c_1\alpha (t - t_0) \right)^{1/\alpha}} = \mathcal{O}\left(t^{-1/\alpha}\right) & \mbox{otherwise}.
  \end{cases}
\end{equation*}
where
\begin{align*}
t_0 \defeq \max \left\{1, \left\lfloor \log_{1 - c_1} \left(\frac{(c_1/c_2)^{1/\alpha}}{c_0} \right)\right\rfloor \right\}.
\end{align*}
\end{lemma}
This allows us to convert the per-iteration contractions to convergence rates.

\begin{theorem}
\label{th:FW_generalizedself-concordant_convergence_str_cvx_case_bounded_gradient}
Suppose $\cx$ is a compact $(\kappa, q)$-uniformly convex set and $f$ is a $\left(M, \nu\right)$ generalized self-concordant function with $\nu \geq 2$.
Furthermore, assume that $\min_{\vx \in \cx} \norm{\nabla f(\vx)} \geq C$.
Then, the Frank-Wolfe algorithm with Backtrack (\cref{fw2}) results in a convergence:
\begin{equation*}
  h_t \leq  \begin{cases}
     h(\vx_0) \left( 1 - \frac{1}{2} \min\left\{ 1,\frac{\kappa C}{\tilde{L}} \right\} \right)^t   & \mbox{if } q = 2 \\
     \frac{h(\vx_0)}{2^t} & \mbox{if } q > 2, 1 \leq t\leq t_0 \\
     \frac{(\tilde{L}^{q} / (\kappa C)^{2})^{1/(q-2)}}{(1 + (q-2)(t - t_0)/(2q))^{q/(q-2)}} = \mathcal{O}\left(t^{-q/(q-2)}\right) & \mbox{if } q > 2, t > t_0,
  \end{cases}
\end{equation*}
for $t\geq 1$, where:
\begin{align*}
t_0 = \max \left\{1, \left\lfloor \log_{1/2}\left( \frac{(\tilde{L}^{q} / (\kappa C)^{2})^{1/(q - 2)}}{h(\vx_0)}\right)\right\rfloor \right\}.
\end{align*}
and $\tilde{L} \defeq \max\{ \tau L_f^{\cl_0}, L_{-1}\}$, where $\tau>1$ is the backtracking parameter, $L_{-1}$ is the initial smoothness estimate in \cref{backtrack}, and
\begin{align*}
  L_f^{\cl_0} = \max\limits_{\vu \in \cl_0, \vd \in \rr^n} \norm{\vd}_{\nabla^2 f(\vu)}^2/\norm{\vd}^2_2.
\end{align*}
\begin{proof}
At iteration $t$, the backtracking line search strategy finds through successive function evaluations a $L_t > 0$ such that:
\begin{align*}
h(\vx_{t+1}) & \leq h(\vx_t) - \gamma_t g(\vx_t) + \frac{L_t \gamma_t^2}{2} \norm{\vx_t - \mathbf{v}_t}^2.
\end{align*}
Finding the $\gamma_t$ that maximizes the right-hand side of the previous inequality leads to:
\begin{align*}
\gamma_t = \min \{1, g(\vx_t)/(L_t \norm{\vx_t - \mathbf{v}_t}^2) \},
\end{align*}
which is the step size ultimately taken by the algorithm at iteration $t$. Note that if $\gamma_t = 1$ this means that $g(\vx_t) \geq L_t \norm{\vx_t - \mathbf{v}_t}^2$, which when plugged into the inequality above leads to $h(\vx_{t+1}) \leq h(\vx_t)/2$. Conversely, for $\gamma_t < 1$ we have that $h(\vx_{t+1}) \leq h(\vx_t) - g(\vx_t)^2/(2L_t \norm{\vx_t - \mathbf{v}_t}^2)$. Focusing on this case and using the bounds $g(\vx_t)\geq h(\vx_t)$ and $g(\vx_t)\geq \kappa \norm{\nabla f(\vx_t)} \norm{\vx_t - \mathbf{v}_t}^q$ from Proposition~\ref{appx:key_proposition_str_cvx_sets} leads to:
\begin{align}
h(\vx_{t+1}) & \leq h(\vx_t)  - h(\vx_t)^{2-2/q} \frac{ \left(\kappa \norm{\nabla f(\vx_t)}\right)^{2/q} }{2L_t}  \label{appx:eq:smoothness_bound} \\
& \leq  h(\vx_t)  - h(\vx_t)^{2-2/q} \frac{ \left(\kappa C\right)^{2/q} }{2\tilde{L}},
\end{align}
where the last inequality simply comes from the bound on the gradient norm, and the fact that $L_t \leq \tilde{L}$, for  $\tilde{L} \defeq \max\{ \tau L_f^{\cl_0}, L_{-1}\}$, where $\tau>1$ is the backtracking parameter and 
$L_{-1}$ is the initial smoothness estimate in \cref{backtrack}. Reordering this expression and putting together the two cases we have that:
\begin{align*}
h(\vx_t) - h(\vx_{t+1}) & \geq h(\vx_t) \min\left\{ \frac{1}{2}, \frac{ \left(\kappa C\right)^{2/q} }{2\tilde{L}} h(\vx_t)^{1 - 2/q} \right\} .
\end{align*}
For the case where $q = 2$ we get a linear contraction in primal gap. Using Lemma~\ref{appx:contraction_to_convergence_rate} to go from a contraction to a convergence rate for $q >2$ we have that:
\begin{equation*}
  h_t \leq  \begin{cases}
     h(\vx_0) \left( 1 - \frac{1}{2} \min\left\{ 1,\frac{\kappa C}{\tilde{L}} \right\} \right)^t   & \mbox{if } q = 2 \\
     \frac{h(\vx_0)}{2^t} & \mbox{if } q > 2, 1 \leq t\leq t_0 \\
     \frac{(\tilde{L}^{q} / (\kappa C)^{2})^{1/(q-2)}}{(1 + (q-2)(t - t_0)/(2q))^{q/(q-2)}} = \mathcal{O}\left(t^{-q/(q-2)}\right) & \mbox{if } q > 2, t > t_0,
  \end{cases}
\end{equation*}
for $t\geq 1$, where:
\begin{align*}
t_0 = \max \left\{1, \left\lfloor \log_{1/2}\left( \frac{(\tilde{L}^{q} / (\kappa C)^{2})^{1/(q - 2)}}{h(\vx_0)}\right)\right\rfloor \right\},
\end{align*}
which completes the proof.
\end{proof}
\end{theorem}

However, in the general case, we cannot assume that the norm of the gradient is bounded away from zero over $\cx$. We deal with the general case using local strong convexity in \cref{th:FW_generalizedself-concordant_convergence_str_cvx_case}.

\begin{theorem}
\label{th:FW_generalizedself-concordant_convergence_str_cvx_case}
Suppose $\cx$ is a compact $(\kappa, q)$-uniformly convex set and $f$ is a $\left(M, \nu\right)$ generalized self-concordant function with $\nu \geq 2$ for which domain does not contain any straight line.
Then, the Frank-Wolfe algorithm with Backtrack (\cref{fw2}) results in a convergence:
\begin{equation*}
  h_t \leq  \begin{cases}
     \frac{h(\vx_0)}{2^t} & \mbox{if } 1 \leq t\leq t_0 \\
     \frac{(\tilde{L}^q /(\kappa^2 \mu_f^{\cl_0}))^{1/(q-1)}}{(1 + (q-1)(t - t_0)/(2q))^{q/(q-1)}} = \mathcal{O}\left(t^{-q/(q-1)}\right) & \mbox{if } t > t_0,
  \end{cases}
\end{equation*}
for $t\geq 1$, where:
\begin{align*}
\cl_0 &= \left\{ \vx \in \dom f \cap \cx \mid f(\vx) \leq f(\vx_0)\right\} \\
t_0 &= \max \left\{1, \left\lfloor \log_{1/2}\left( \frac{(\tilde{L}^q /(\kappa^2 \mu_f^{\cl_0}))^{1/(q-1)}}{h(\vx_0)}\right)\right\rfloor \right\}
\end{align*}
and $\tilde{L} \defeq \max\{ \tau L_f^{\cl_0}, L_{-1}\}$, where $\tau>1$ is the backtracking parameter, $L_{-1}$ is the initial smoothness estimate in \cref{backtrack}.


\begin{proof}
As the algorithm makes monotonic primal progress we have that $\vx_t \in \cl_0$ for $t\geq 0$.
The proof proceeds very similarly as before, except for the fact that now we have to bound $\norm{\nabla f(\vx_t)}$ using $\mu_f^{\cl_0}$-strong convexity for points $\vx_t, \vx_t + \gamma_t (\mathbf{v}_t - \vx_t) \in \cl_0$.
Continuing from \eqref{appx:eq:smoothness_bound} for the case where $\gamma_t < 1$ and using the fact that, given $f^{u}$ the unconstrained minimum of $f$, $h(\vx_t) \leq f(\vx_t) - f^{u} \leq \norm{\nabla f(\vx_t)}^2/(2\mu_f^{\cl_0})$ we have that:
\begin{align*}
h(\vx_{t+1}) & \leq h(\vx_t)  - h(\vx_t)^{2-2/q} \frac{ \left(\kappa \norm{\nabla f(\vx_t)}\right)^{2/q} }{2L_t} \\
& \leq h(\vx_t)  - h(\vx_t)^{2-1/q} \frac{ \kappa^{2/q} (\mu_f^{\cl_0})^{1/q} 2^{1/q - 1} }{\tilde{L}},
\end{align*}
where we have also used the bound $L_t \leq \tilde{L}$ in the last equation. This leads us to a contraction, together with the case where $\gamma_t = 1$, which is unchanged from the previous proofs, of the form:
\begin{align*}
h(\vx_t) - h(\vx_{t+1}) & \geq h(\vx_t) \min\left\{ \frac{1}{2}, \frac{ \kappa^{2/q}  (\mu_f^{\cl_0})^{1/q} 2^{1/q - 1} }{\tilde{L}} h(\vx_t)^{1 - 1/q} \right\}.
\end{align*}
Using again Lemma~\ref{appx:contraction_to_convergence_rate} to go from a contraction to a convergence rate for $q >2$ we have that:
\begin{equation*}
  h_t \leq  \begin{cases}
     \frac{h(\vx_0)}{2^t} & \mbox{if } 1 \leq t\leq t_0 \\
     \frac{(\tilde{L}^q /(\kappa^2 \mu_f^{\cl_0}))^{1/(q-1)}}{(1 + (q-1)(t - t_0)/(2q))^{q/(q-1)}} = \mathcal{O}\left(t^{-q/(q-1)}\right) & \mbox{if } t > t_0,
  \end{cases}
\end{equation*}
for $t\geq 1$, where:
\begin{align*}
t_0 = \max \left\{1, \left\lfloor \log_{1/2}\left( \frac{(\tilde{L}^q /(\kappa^2 \mu_f^{\cl_0}))^{1/(q-1)}}{h(\vx_0)}\right)\right\rfloor \right\},
\end{align*}
which completes the proof.
\end{proof}
\end{theorem}

In Table~\ref{table:bound_comparison_backtracking_versions} we provide an oracle complexity breakdown for the Frank-Wolfe algorithm with Backtrack (\texttt{B-FW}), also referred to as \texttt{LBTFW-GSC} in \cite{dvurechensky2020generalized}, when minimizing over a $(\kappa, q)$-uniformly convex set. 

\begin{table}[H]\centering
  \footnotesize
\begin{tabular}{lcccc}
\toprule
\textbf{Algorithm}  &  \textbf{Assumptions} & \textbf{Oracle calls}  & \textbf{Reference}\\
\midrule
\texttt{B-FW}/\texttt{LBTFW-GSC}$^{\ddagger}$  &  $\vx^* \in \interior \left(\cx \cap \dom f \right) $  & $\mathcal{O}(\log 1/\varepsilon)$ &  \textbf{This work} \\
\texttt{B-FW}/\texttt{LBTFW-GSC}$^{\ddagger}$ &  $\min_{\vx \in \cx} \norm{\nabla f(\vx)} > 0$, $q = 2$  & $\mathcal{O}(\log 1/\varepsilon)$ &  \textbf{This work} \\
\texttt{B-FW}/\texttt{LBTFW-GSC}$^{\ddagger}$ &  $\min_{\vx \in \cx} \norm{\nabla f(\vx)} > 0$, $q > 2$ & $\mathcal{O}\left(\varepsilon^{-(q-2)/q}\right)$ &  \textbf{This work}  \\
\texttt{B-FW}/\texttt{LBTFW-GSC}$^{\ddagger}$ &  No straight lines in $\dom{f}$ & $\mathcal{O}\left(\varepsilon^{-(q-1)/q}\right)$ &  \textbf{This work}  \\
\bottomrule
\end{tabular}
\smallskip
\caption{
\textbf{Complexity comparison for \texttt{B-FW} (\cref{fw2}) when minimizing over a $(\kappa, q)$-uniformly convex set}: Number of iterations needed to reach an $\varepsilon$-optimal solution in $h(\vx)$ for Problem~\ref{prob:main} in several cases of interest. We use the superscript $\ddagger$ to indicate that constants in the convergence bounds depend on user-defined inputs. Oracle calls refer simultaneously to FOO, ZOO, LMO, and DO calls.}
 \label{table:bound_comparison_backtracking_versions}
\end{table}

\section{Away-step and Blended Pairwise Conditional Gradients}

When the domain $\cx$ is a polytope, one can obtain linear convergence in primal gap for a generalized self-concordant function using the well known \emph{Away-step Frank-Wolfe} (AFW) algorithm \citep{guelat1986some,lacoste15} shown in \cref{appx:afw}
and the more recent \emph{Blended Pairwise Conditional Gradients} (BPCG) algorithm \citep{tsuji2022pairwise} with the adaptive step size of \citet{pedregosa2018step}.
We use $\cs_t$ to denote the \emph{active set} at iteration $t$, that is, the set of vertices of the polytope that gives rise to $\vx_t$ as a convex combination with positive weights.

For AFW, we can see that the algorithm either chooses to perform what is known as a \emph{Frank-Wolfe} step in Line~\ref{alg:appx:fw_step} of \cref{appx:afw}
if the Frank-Wolfe gap $g(\vx)$ is greater than the \emph{away gap} $\innp{\nabla f(\vx_t),\va_t - \vx_t}$ or an \emph{Away} step in \cref{alg:appx:afw_step} of \cref{appx:afw} otherwise.
Similarly for BPCG, the algorithm performs a Frank-Wolfe step in Line~\ref{alg:BPCG:if_condition} of \cref{alg:BPCG} if the Frank-Wolfe gap is greater than the pairwise gap $\innp{\nabla f(\vx_t),\va_t - \vs_t}$.
For simplicity of exposition, we made both algorithms start with a vertex of of $\cx$ in $\dom{f}$.
Although this could be too restrictive for some applications (e.g.~when the function includes a barrier of the polytope), it is easy to warm-start the active set with an initial combination of vertices.

\begin{algorithm}[H]
	\caption{Away-step Frank-Wolfe (\texttt{B-AFW}) with the step size of \citet{pedregosa2018step}}
\label{appx:afw}  
\footnotesize
\begin{algorithmic}[1]
\Require Vertex $\vx_0 \in \dom{f}$ of $\cx$, function \(f\), initial smoothness estimate $L_{-1}$
\hrulealg
  \State $\cs_0 \leftarrow \{\vx_0\}$, $\vlambda_0 \leftarrow \{1\}$
\For{$t=0$ \textbf{to} $\dotsc$}
  \State $\mathbf{v}_t \leftarrow \argmin_{\mathbf{v} \in  \mathcal{X}} \innp{\nabla f\left(\vx_t \right), \mathbf{v}}$
  \State $\va_t \leftarrow \argmax_{\mathbf{v} \in  \cs_t} \innp{\nabla f\left(\vx_t \right), \mathbf{v}}$
   \If{$\innp{\nabla f(\vx_t),\vx_t - \mathbf{v}_t} \geq \innp{\nabla f(\vx_t),\va_t - \vx_t}$} \label{alg:appx:if_condition}
   \State $\vd_t \leftarrow \mathbf{v}_t - \vx_t$, $\gamma_{\max} \leftarrow 1$ \label{alg:appx:fw_step}
   \Else
   \State $\vd_t \leftarrow \vx_t - \va_t$, $\gamma_{\max} \leftarrow \vlambda_t(\va_t)/\left( 1 - \vlambda_t(\va_t)\right)$\label{alg:appx:afw_step}
   \EndIf
     \State $\gamma_t, L_{t} \leftarrow$ \texttt{Backtrack}$\left(f, \vx_t, \vd_t, \nabla f(\vx_t), L_{t-1}, \gamma_{\max} \right)$
   \State $\vx_{t+1} \leftarrow \vx_t + \gamma_t \vd_t$
   \State Update $\cs_t$ and $\vlambda_t$ to $\cs_{t+1}$ and $\vlambda_{t+1}$
\EndFor
\end{algorithmic}
\end{algorithm}

\begin{algorithm}[H]
	\caption{Blended Pairwise Conditional Gradients (\texttt{B-BPCG}) with the step size of \citet{pedregosa2018step}}\label{alg:BPCG}
  \footnotesize
\begin{algorithmic}[1]
\Require Vertex $\vx_0 \in \dom{f}$ of $\cx$, function \(f\), initial smoothness estimate $L_{-1}$
\hrulealg
  \State $\cs_0 \leftarrow \{\vx_0\}$, $\vlambda_0 \leftarrow \{1\}$
\For{$t=0$ \textbf{to} $\dotsc$}
  \State $\mathbf{v}_t \leftarrow \argmin_{\mathbf{v} \in  \mathcal{X}} \innp{\nabla f\left(\vx_t \right), \mathbf{v}}$
  \State $\va_t \leftarrow \argmax_{\mathbf{v} \in  \cs_t} \innp{\nabla f\left(\vx_t \right), \mathbf{v}}$
  \State $\mathbf{s}_t \leftarrow \argmin_{\mathbf{v} \in  \cs_t} \innp{\nabla f\left(\vx_t \right), \mathbf{v}}$
   \If{$\innp{\nabla f(\vx_t),\vx_t - \mathbf{v}_t} \geq \innp{\nabla f(\va_t),\mathbf{s}_t - \va_t}$} \label{alg:BPCG:if_condition}
   \State $\vd_t \leftarrow \mathbf{v}_t - \vx_t$, $\gamma_{\max} \leftarrow 1$\label{alg:BPCG:fw_step}
   \Else
   \State $\vd_t \leftarrow \va_t - \mathbf{s}_t$, $\gamma_{\max} \leftarrow \vlambda_t(\va_t)$ \label{alg:appx:pw_step}
   \EndIf
     \State $\gamma_t, L_{t} \leftarrow$ \texttt{Backtrack}$\left(f, \vx_t, \vd_t, \nabla f(\vx_t), L_{t-1}, \gamma_{\max} \right)$
   \State $\vx_{t+1} \leftarrow \vx_t + \gamma_t \vd_t$
   \State Update $\cs_t$ and $\vlambda_t$ to $\cs_{t+1}$ and $\vlambda_{t+1}$
\EndFor
\end{algorithmic}
\end{algorithm}

Both proofs of linear convergence follow closely from \cite{pedregosa2018step,lacoste15} but leverage generalized self-concordant instead of smoothness and strong convexity.
One of the key inequalities used in the proof is a scaling inequality from \citet{lacoste15} very similar to the one shown in Proposition~\ref{th:bound_interior_optimum} and Proposition~\ref{appx:key_proposition_str_cvx_sets}, which we state next:

\begin{proposition}\label{appx:key_proposition_polytope}
Let $\cx\subseteq \rr^n$ be a polytope, and denote by $\cs$ the set of vertices of the polytope $\cx$ that gives rise to $\vx\in \cx$ as a convex combination with positive weights, then for all $\vy \in \cx$:
\begin{align*}
\innp{\nabla f(\vx), \va - \mathbf{v}} \geq \delta \frac{\innp{\nabla f(\vx), \vx - \vy}}{\norm{ \vx - \vy}},
\end{align*}
where $\mathbf{v} = \argmin_{\vu\in \cx}\innp{\nabla f(\vx), \vu}$, $\va = \argmax_{\vu\in \cs}\innp{\nabla f(\vx), \vu}$, and $\delta >0$ is the \emph{pyramidal width} of $\cx$.
\end{proposition}

\begin{theorem}
\label{appx:th:AFW_generalizedself-concordant_convergence_polytope}
Suppose $\cx$ is a polytope and $f$ is a $\left(M, \nu\right)$ generalized self-concordant function with $\nu \geq 2$ for which the domain does not contain any straight line.
Then, both AFW and BPCG with Backtrack achieve a convergence rate:
\begin{align*}
h(\vx_{t}) \leq h(\vx_0)\left( 1 -  \frac{\mu_f^{\cl_0}}{4\tilde{L}}  \left(\frac{\delta}{D}\right)^2\right)^{\left\lceil (t-1)/2\right\rceil},
\end{align*}
where $\delta$ is the pyramidal width of the polytope $\cx$, $\tilde{L} \defeq \max\{ \tau L_f^{\cl_0}, L_{-1}\}$, $\tau>1$ is the backtracking parameter, 
$L_{-1}$ is the initial smoothness estimate in \cref{backtrack}.
\begin{proof}
Proceeding very similarly as in the proof of Theorem~\ref{th:FW_generalizedself-concordant_convergence_str_cvx_case_inner_opt},
we have that as the backtracking line search makes monotonic primal progress, we know that for $t\geq 0$ we will have that $\vx_t \in \cl_0$.
As the function is $\mu_f^{\cl_0}$-strongly convex over $\cl_0$, we can use the appropriate inequalities from strong convexity in the progress bounds.
Using this aforementioned property, together with the scaling inequality of Proposition~\ref{appx:key_proposition_polytope} results in:
\begin{align}
h(\vx_t) = f(\vx_t) - f(\vx^*) & \leq \frac{\innp{\nabla f(\vx_t), \vx_t - \vx^*}^2}{2\mu_f^{\cl_0}\norm{\vx_t - \vx^*}^2}\nonumber \\
& \leq \frac{\innp{\nabla f(\vx_t), \va_t - \mathbf{v}_t}^2}{2\mu_f^{\cl_0}\delta^2}.\label{eq:firstbound}
\end{align}
The first inequality comes from the $\mu_f^{\cl_0}$-strong convexity over $\cl_0$ (see, e.g., \cite[Lemma 2.13]{cg_survey}), and the second inequality comes from applying
Proposition~\ref{appx:key_proposition_polytope} with $\vy = \vx^*$.
For AFW, we can expand the expression of the numerator of the bound in \eqref{eq:firstbound}:
\begin{align}
  f(\vx_t) - f(\vx^*) & \leq \frac{\left(\innp{\nabla f(\vx_t), \va_t - \vx_t} + \innp{\nabla f(\vx_t), \vx_t - \mathbf{v}_t}\right)^2}{2\mu_f^{\cl_0}\delta^2}. \label{eq:key:inequality:afw}
\end{align}
Note that if the Frank-Wolfe step is chosen in Line~\ref{alg:appx:fw_step}, then:
\begin{align*}
  -\innp{\nabla f(\vx_t),\vd_t} = \innp{\nabla f(\vx_t),\vx_t - \mathbf{v}_t} \geq \innp{\nabla f(\vx_t),\va_t - \vx_t},  
\end{align*}
otherwise, if an away step is chosen in Line~\ref{alg:appx:afw_step}, then:
\begin{align*}
  \innp{\nabla f(\vx_t),\vx_t - \mathbf{v}_t} < \innp{\nabla f(\vx_t),\va_t - \vx_t} = -\innp{\nabla f(\vx_t),\vd_t}.  
\end{align*}
In both cases, we have that:
\begin{align}
f(\vx_t) - f(\vx^*) & \leq \frac{2\innp{\nabla f(\vx_t), \vd_t}^2}{\mu_f^{\cl_0}\delta^2}. \label{appx:key:inequality1}
\end{align}
For BPCG, we can directly exploit \cite[Lemma 3.5]{tsuji2022pairwise}, which establishes the following bound at every iteration of the algorithm:
\begin{align*}
  -2\innp{\nabla f(\vx_t), \vd_t} \geq \innp{\nabla f(\vx_t), \va_t - \mathbf{v}_t},
\end{align*}
resulting in the same inequality \eqref{appx:key:inequality1}.\footnote{Note the minus sign accounting for the difference in the definition of $\vd_t$ between our paper and \citet{tsuji2022pairwise}.}
Note that using a similar reasoning, as $\innp{\nabla f(\vx_t),\vx_t - \mathbf{v}_t} = g(\vx_t)$, in both cases it holds that:
\begin{align}
h(\vx_t)\leq g(\vx_t) \leq -\innp{\nabla f(\vx_t),\vd_t}. \label{appx:key:inequality2}
\end{align}
As in the preceding proofs, the backtracking line search in \cref{backtrack} will either output a point $\gamma_t = \gamma_{\max}$ or $\gamma_t < \gamma_{\max}$.
In any case, for both AFW and BPCG, and regardless of the type of step taken, \cref{backtrack} will find and output a smoothness estimate $L_t$ and a step size $\gamma_t$ such that:
\begin{align}
h(\vx_{t+1}) - h(\vx_t) \leq  \frac{L_t\gamma_t^2}{2}  \norm{\vd_t}^2 + \gamma_t \innp{\nabla f(\vx_t), \vd_t}. \label{appx:local:smoothness:AFW}
\end{align}
As before, we will have two cases differentiating whether the step size $\gamma_t$ is maximal.
If $\gamma_t = \gamma_{\max}$ we know by observing Line~\ref{backtrack:stepsize} of \cref{backtrack} that:
\begin{align*}
  -\innp{\nabla f(\vx_t), \vd_t} \geq \gamma_{\max}L_t \norm{\vd_t}^2,  
\end{align*}
which combined with \eqref{appx:local:smoothness:AFW} results in:
\begin{align*}
  h(\vx_{t+1}) - h(\vx_{t}) \leq \innp{\nabla f(\vx_t), \vd_t}\frac{\gamma_{\max}}{2}.
\end{align*}
In the case where $\gamma_t < \gamma_{\max}$, we have:
\begin{align*}
  & -\innp{\nabla f(\vx_t), \vd_t} < \gamma_{\max} L_t \norm{\vd_t}^2,\\
  & \gamma_t = -\innp{\nabla f(\vx_t), \vd_t}/(L_t\norm{\vd_t}^2).
\end{align*}
Plugging the expression of $\gamma_t$ into \eqref{appx:local:smoothness:AFW} yields
\begin{align*}
  h(\vx_{t+1}) - h(\vx_{t}) \leq -\innp{\nabla f(\vx_t), \vd_t}^2/(2L_t\norm{\vd_t}^2).  
\end{align*}
In any case, we can rewrite \eqref{appx:local:smoothness:AFW} as:
\begin{align}
h(\vx_t) - h(\vx_{t+1}) &\geq   \min\left\{\frac{-\innp{\nabla f(\vx_t), \vd_t}\gamma_{\max}}{2}, \frac{\innp{\nabla f(\vx_t), \vd_t}^2}{2L_t\norm{\vd_t}^2} \right\}. \label{appx:lin:conv:ineq}
\end{align}
We can now use the inequality in \eqref{appx:key:inequality1} to bound the second term in the minimization component of \eqref{appx:lin:conv:ineq}, and \eqref{appx:key:inequality2} to bound the first term. This leads to:
\begin{align}
h(\vx_t) - h(\vx_{t+1}) &\geq   h(\vx_t)\min\left\{\frac{\gamma_{\max}}{2}, \frac{\mu_f^{\cl_0}\delta^2}{4L_t\norm{\vd_t}^2} \right\} \\
&\geq   h(\vx_t)\min\left\{\frac{\gamma_{\max}}{2}, \frac{\mu_f^{\cl_0}\delta^2}{4\tilde{L}D^2} \right\},
\end{align}
where in the last inequality we use $\norm{\vd_t}\leq D$ and $L_t \leq \tilde{L}$ for all $t$.
It remains to bound $\gamma_{\max}$ away from zero to obtain the linear convergence bound.
For Frank-Wolfe steps, we immediately have $\gamma_{\max} = 1$, but for away or pairwise steps, there is no straightforward way of bounding $\gamma_{\max}$ away from zero.
One of the key insights from \citet{lacoste15} is that instead of bounding $\gamma_{\max}$ away from zero for all steps up to iteration $t$, we can instead bound the number of away steps with a step size $\gamma_t = \gamma_{\max}$ up to iteration $t$, which are steps that reduce the cardinality of the active set $\cs_t$ and satisfy $h(\vx_{t+1}) \leq h(\vx_{t})$.
The same argument is used in \citet{tsuji2022pairwise} to prove the convergence of BPCG.
This leads us to consider only the progress provided by the remaining steps, which are Frank-Wolfe steps and away steps for AFW or pairwise steps for BPCG with $\gamma_t < \gamma_{\max}$.
For a number of steps $t$, only at most half of these steps could have been away steps with $\gamma_t = \gamma_{\max}$, as we cannot drop more vertices from the active set than the number of vertices we could have potentially picked up with Frank-Wolfe steps.
For the remaining $\left\lceil (t-1)/2\right\rceil$ steps, we know that:
\begin{align*}
  h(\vx_t) - h(\vx_{t+1})\geq h(\vx_t) \frac{\mu_f^{\cl_0}\delta^2}{4\tilde{L}D^2}.
\end{align*}
Therefore, we have that the primal gap satisfies:
\begin{align*}
h(\vx_{t}) \leq h(\vx_0)\left( 1 -  \frac{\mu_f^{\cl_0}\delta^2}{4\tilde{L}D^2}\right)^{\left\lceil (t-1)/2\right\rceil}.
\end{align*}
\end{proof}
\end{theorem}

We can make use of the proof of convergence in primal gap to prove linear convergence in Frank-Wolfe gap. In order to do so, we recall a quantity formally defined in \cite{kerdreux2019restarting} but already implicitly used earlier in \cite{lacoste15} as:
\begin{align*}
w(\vx_t, \cs_t) & \defeq \max\limits_{\vu \in \cs_t, \mathbf{v} \in \cx} \innp{\nabla f(\vx_t), \vu - \mathbf{v}} \\
& = \max\limits_{\vu \in \cs_t} \innp{\nabla f(\vx_t), \vu - \vx_t} + \max\limits_{ \mathbf{v} \in \cx} \innp{\nabla f(\vx_t), \vx_t - \mathbf{v}} \\
& = \max\limits_{\vu \in \cs_t} \innp{\nabla f(\vx_t), \vu - \vx_t} + g(\vx_t).
\end{align*}
Note that $w(\vx_t, \cs_t)$ provides an upper bound on the Frank-Wolfe gap as the first term in 
the definition, the so-called \emph{away gap}, is positive.

\begin{theorem}
\label{appx:th:AFW_generalizedself-concordant_convergence_polytope_FW_gap}
Suppose $\cx$ is a polytope and $f$ is a $\left(M, \nu\right)$ generalized self-concordant function with $\nu \geq 2$ for which the domain does not contain any straight line.
Then, AFW and BPCG with Backtrack both contract the Frank-Wolfe gap linearly, i.e., $\min_{1 \leq t \leq T} g(\vx_t) \leq \varepsilon$ after $T = \mathcal O(\log 1/\varepsilon)$ iterations.
\begin{proof}
We observed in the proof of \cref{appx:th:AFW_generalizedself-concordant_convergence_polytope} that regardless of the type of step chosen in AFW and BPCG, the following holds:
\begin{align*}
  -2\innp{\nabla f(\vx_t), \vd_t} \geq  \innp{\nabla f(\vx_t), \vx_t - \mathbf{v}_t} +  \innp{\nabla f(\vx_t), \va_t - \vx_t} = w(\vx_t, \cs_t).  
\end{align*}
On the other hand, we also have that $h(\vx_t) - h(\vx_{t+1})\leq h(\vx_t)$. Plugging these bounds into the right-hand side and the left hand side of \ref{appx:lin:conv:ineq} in Theorem~\ref{appx:th:AFW_generalizedself-concordant_convergence_polytope}, and using the fact that $\norm{\vd_t} \leq D$ we have that:
\begin{align*}
\min\left\{\frac{w(\vx_t, \cs_t)\gamma_{\max}}{4}, \frac{w(\vx_t, \cs_t)^2}{8L_t D^2} \right\}  \leq h(\vx_{t})   \leq h(\vx_0)\left( 1 -  \frac{\mu_f^{\cl_0}}{4\tilde{L}}  \left(\frac{\delta}{D}\right)^2\right)^{\left\lceil (t-1)/2\right\rceil},
\end{align*}
where the second inequality follows from the convergence bound on the primal gap from Theorem~\ref{appx:th:AFW_generalizedself-concordant_convergence_polytope}. Considering the steps that are not away steps with $\gamma_t = \gamma_{\max}$ as in the proof of Theorem~\ref{appx:th:AFW_generalizedself-concordant_convergence_polytope}, leads us to:
\begin{align*}
	g(\vx_t) \leq w(\vx_t, \cs_t) \leq 4 h(\vx_0) \max\left\{1 , \sqrt{ \frac{\tilde{L}D^2}{2h(\vx_0) }  } \right\} \left( 1 -  \frac{\mu_f^{\cl_0}}{4\tilde{L}}  \left(\frac{\delta}{D}\right)^2\right)^{\left\lfloor (t-1)/4\right\rfloor}.
\end{align*}
\end{proof}
\end{theorem}

In Table~\ref{table:bound_comparison_improved} we provide a detailed complexity comparison between the Backtracking AFW (\texttt{B-AFW}) \cref{appx:afw}, and other comparable algorithms in the literature.

\begin{table}[H]\centering
  \footnotesize
\begin{tabular}{lccccc}
\toprule
\textbf{Algorithm} & \textbf{SOO calls} & \textbf{FOO calls} & \textbf{ZOO calls}& \textbf{LMO calls} & \textbf{DO calls} \\
\midrule
\texttt{FW-LLOO} \cite[Alg.7]{dvurechensky2020generalized} & $\mathcal{O}(\log 1/\varepsilon)$ & $\mathcal{O}(\log 1/\varepsilon)$ &  & $\mathcal{O}(\log 1/\varepsilon)$* &  \\
\texttt{ASFW-GSC} \cite[Alg.8]{dvurechensky2020generalized} & $\mathcal{O}(\log 1/\varepsilon)$ & $\mathcal{O}(\log 1/\varepsilon)$ &   & $\mathcal{O}(\log 1/\varepsilon)$ &  \\
\texttt{B-AFW/B-BPCG}$^{\dagger \ddagger}$ &  & $\mathcal{O}(\log 1/\varepsilon)$ &  $\mathcal{O}(\log 1/\varepsilon)$ & $\mathcal{O}(\log 1/\varepsilon)$ & $\mathcal{O}(\log 1/\varepsilon)$ \\
\bottomrule
\end{tabular}
\smallskip
\caption{
\textbf{Complexity comparison}: Number of iterations needed to reach a solution with $h(\vx)$ below $\varepsilon$ for Problem~\ref{prob:main} for Frank-Wolfe-type algorithms in the literature. The asterisk on \texttt{FW-LLOO} highlights the fact that the procedure is different from the standard \texttt{LMO} procedure. The complexity shown for the \texttt{FW-LLOO}, \texttt{ASFW-GSC}, and \texttt{B-AFW} algorithms only apply to polyhedral domains, with the additional requirement that for the former two we need an explicit polyhedral representation of the domain (see Assumption 3 in \cite{dvurechensky2020generalized}), whereas the latter only requires an LMO.
The requirement that we have an explicit polyhedral representation may be limiting, for instance for the matching polytope over non-bipartite graphs, as the size of the polyhedral representation in this case depends exponentially on the number of nodes of the graph \citep{rothvoss2017matching}. We use the superscript $\dagger$ to indicate that the same complexities hold when reaching an $\varepsilon$-optimal solution in $g(\vx)$, and the superscript $\ddagger$ to indicate that constants in the convergence bounds depend on user-defined inputs.}
 \label{table:bound_comparison_improved}
\end{table}

\section{Computational experiments} \label{experiments}

We showcase the performance of the \texttt{M-FW} algorithm,
the second-order step size and the LLOO algorithm from \citet{dvurechensky2020generalized} (denoted by \texttt{GSC-FW} and \texttt{LLOO} in the figures) and the Frank-Wolfe and the Away-step Frank-Wolfe algorithm with the backtracking stepsize of \cite{pedregosa2018step},
denoted by \texttt{B-FW} and \texttt{B-AFW} respectively.
We ran all experiments on a server with 8 Intel Xeon 3.50GHz CPUs
and 32GB RAM in single-threaded mode in Julia 1.6.0 with the \texttt{FrankWolfe.jl} package \citep{besanccon2021frankwolfe}.
The data sets used in the problem instances can be found in \citet{carderera_alejandro_2021_4836009}, the code used for the experiments can be found in \cite{coderepo}.
When running the adaptive step size from \cite{pedregosa2018step},
the only parameter that we need to set is the initial smoothness estimate $L_{-1}$.
We use the initialization proposed in \cite{pedregosa2018step}, namely:
\begin{equation*}
  L_{-1} = \norm{\nabla f(\vx_0) - \nabla f(\vx_0 + \varepsilon (\mathbf{v}_0 - \vx_0))} / (\varepsilon \norm{\mathbf{v}_0 - \vx_0})
\end{equation*}
with $\varepsilon$ set to $10^{-3}$.
The scaling parameters $\tau = 2, \eta = 0.9$ are left at their default values as proposed in \cite{pedregosa2018step} and also used in  \cite{dvurechensky2020generalized}.

We also use the vanilla FW algorithm denoted by \texttt{FW}, which is simply \cref{fw} without Lines~\ref{alg:domain:if} and \ref{alg:same_point} using the traditional $\gamma_t = 2/(t+2)$ open-loop step size rule.
Note that there are no formal convergence guarantees for this algorithm when applied to Problem~\eqref{prob:main}. All figures show the evolution of the $h(\vx_t)$ and $g (\vx_t)$ against $t$ and time with a log-log scale. As in \citet{dvurechensky2020generalized} we implemented the LLOO based variant only for the portfolio optimization instance $\Delta_n$; for the other examples, the oracle implementation was not implemented due to the need to estimate non-trivial parameters. 

As can be seen in all experiments, the Monotonic Frank-Wolfe algorithm is very competitive, outperforming previously proposed variants in both in progress per iteration and time.
The only other algorithm that is sometimes faster is the Away-step Frank-Wolfe variant, which however depends on an active set and can therefore induce up a quadratic time and memory overhead, potentially rendering the method inattractive for very large-scale settings.

\begin{figure}
  \centering
  \begin{tabular}{@{}c@{}}
    \includegraphics[width= 0.45\textwidth]{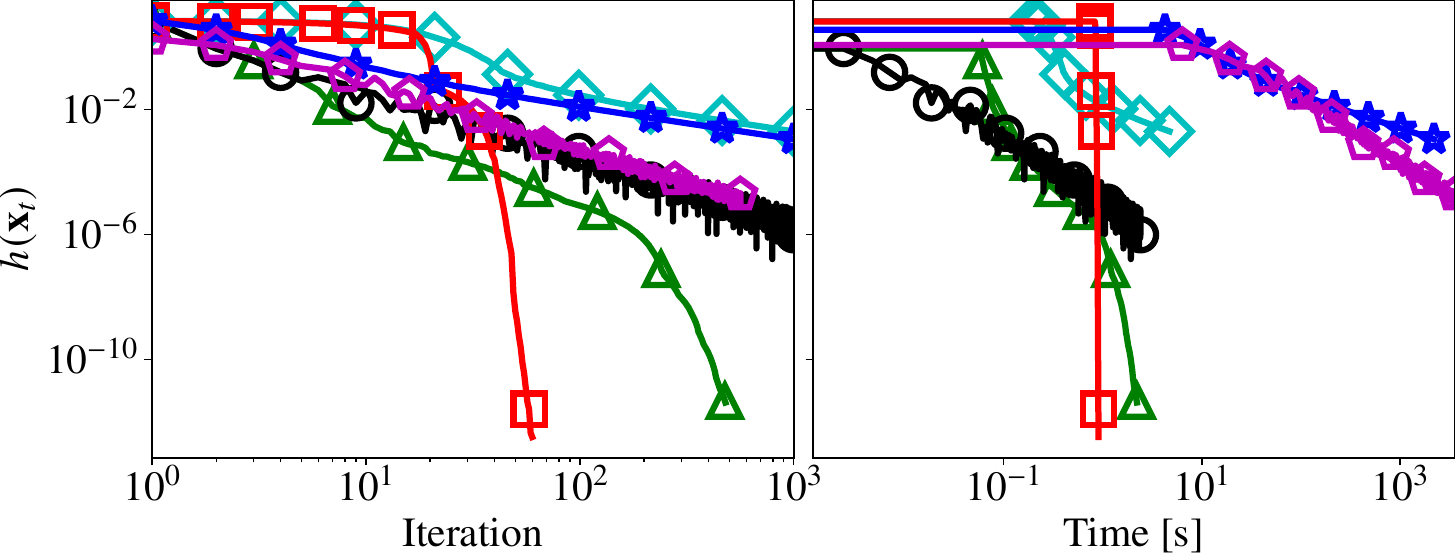} \\[\abovecaptionskip]
  \end{tabular}
  \begin{tabular}{@{}c@{}}
    \includegraphics[width= 0.45\textwidth]{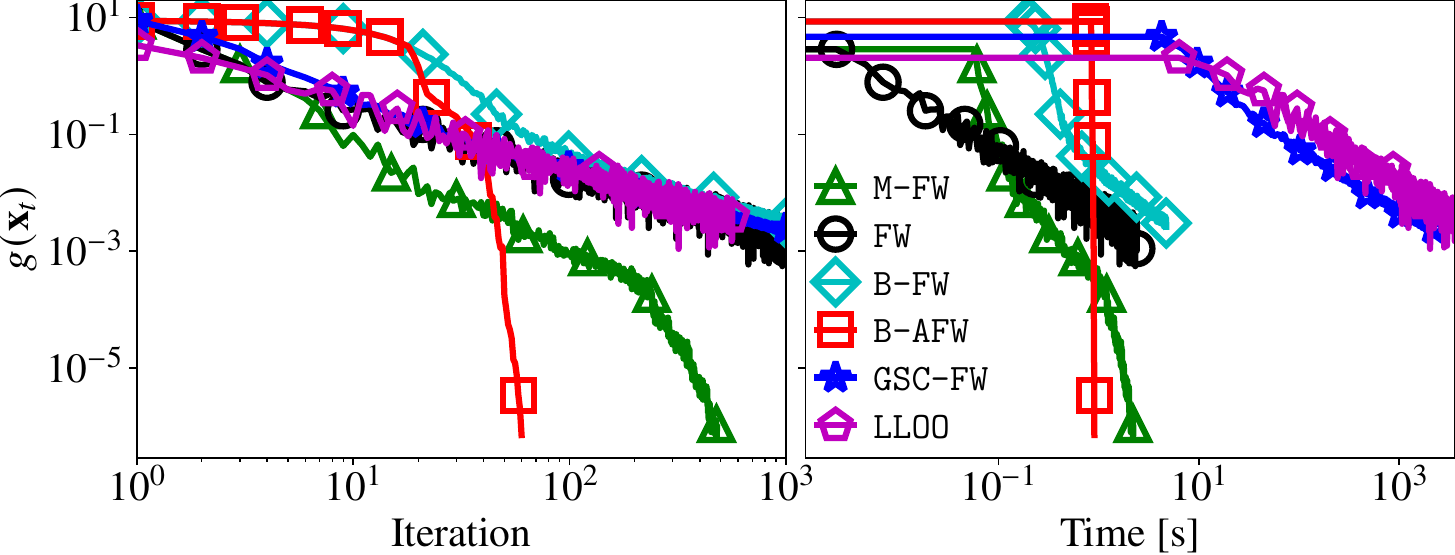} \\[\abovecaptionskip]
  \end{tabular}
  \vspace{-0.25cm}
  \caption{\textbf{Portfolio optimization}: Convergence of $h(\vx_t)$ and $g(\vx_t)$ vs. $t$ and wall-clock time. $n=1000$.}\label{fig:portfolio}
\end{figure}

\begin{figure}
\centering
\begin{subfigure}[b]{0.48\textwidth}
    \centering
    \includegraphics[width= 0.8\textwidth]{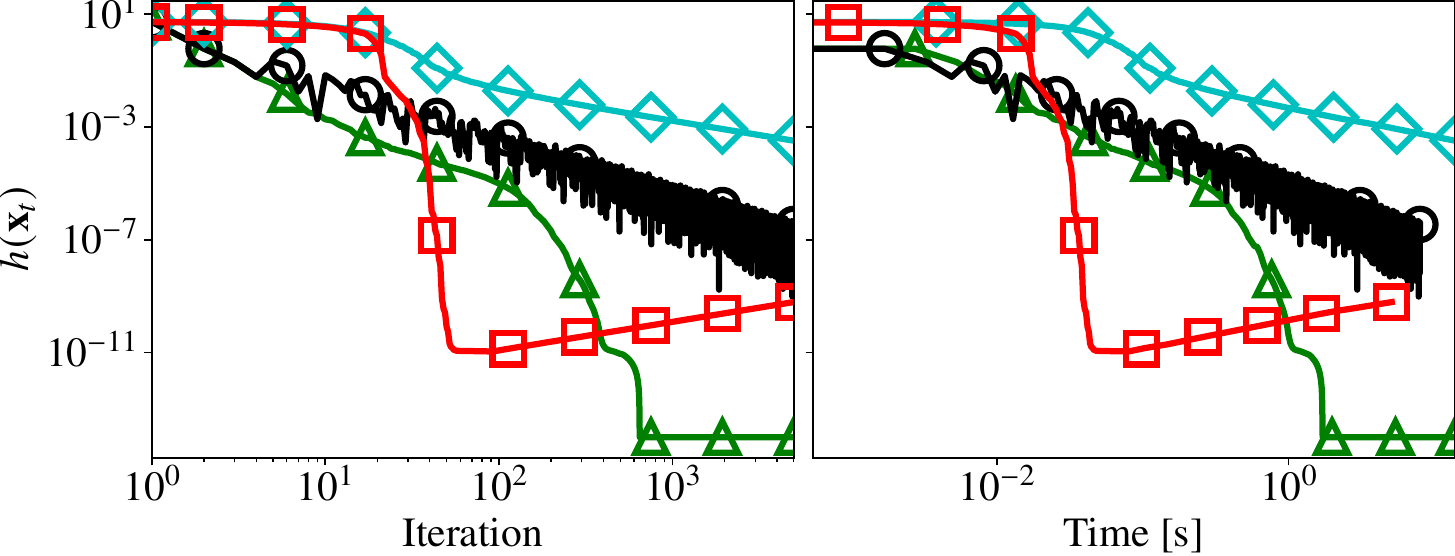}
    \includegraphics[width= 0.8\textwidth]{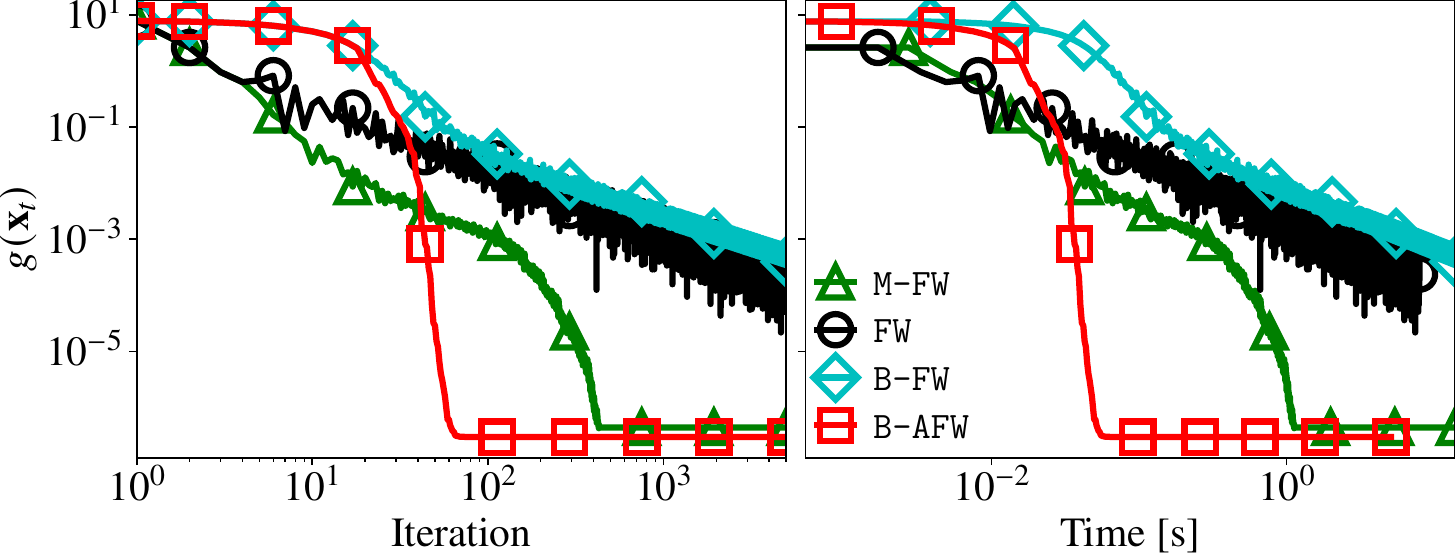}
      \caption{$n=2000$}
\end{subfigure}
\hfill
\begin{subfigure}[b]{0.48\textwidth}
  \centering
  \includegraphics[width= 0.8\textwidth]{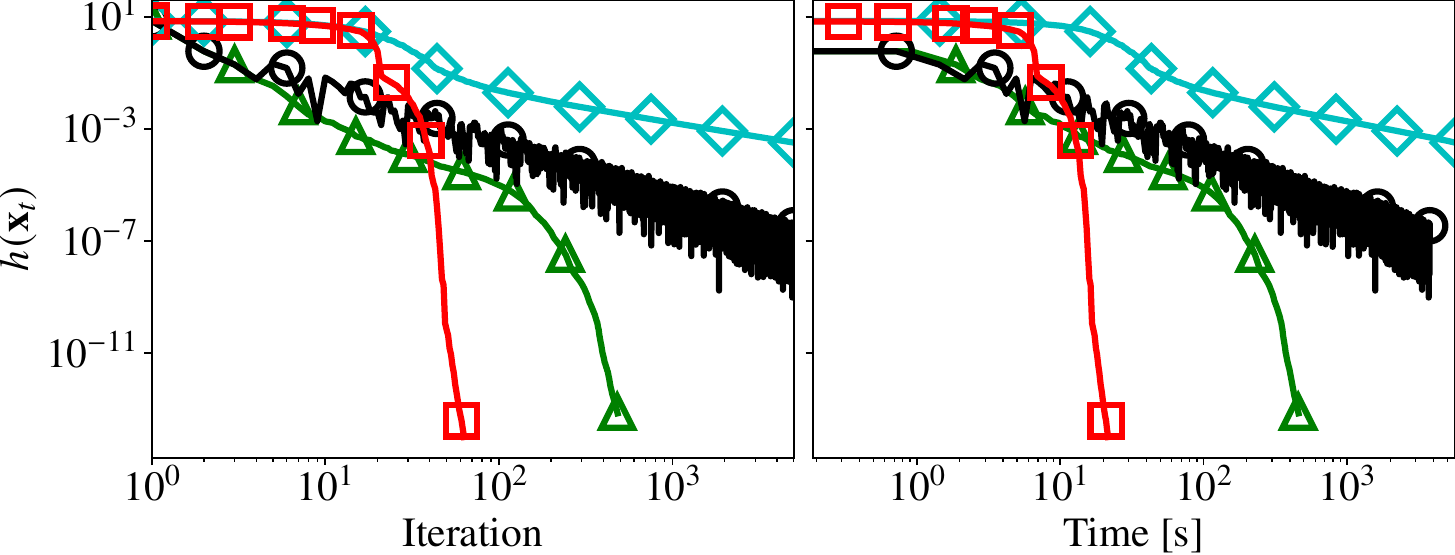}
  \includegraphics[width= 0.8\textwidth]{Images/Primal_gap_Portfolio_appx1_big.pdf}
  \caption{$n=5000$}
  \label{fig:2}
\end{subfigure}
\caption{\textbf{Portfolio optimization}: Convergence of $h(\vx_t)$ and $g(\vx_t)$ vs. $t$ and wall-clock time.}
  \label{fig:portfoliogathered}
\end{figure}

\smallskip\noindent\textbf{Portfolio optimization.}
We consider the portfolio problem with logarithmic returns $f(\vx) = -\sum_{t=1}^p \log(\langle \vr_t, \vx\rangle)$, where
$p$ denotes the number of periods and $\cx = \Delta_n$.
The results are shown in Figure~\ref{fig:portfolio} with all methods and in Figure~\ref{fig:portfoliogathered} on larger instances with first-order methods only.
We use the revenue data $\vr_t$ from \cite{dvurechensky2020generalized}
and add instances generated in a similar fashion from independent Normal random entries with dimension 
1000, 2000, and 5000, and from a Log-normal distribution with $(\mu=0.0, \sigma=0.5)$.  

\begin{figure}
\centering
\begin{subfigure}[b]{0.48\textwidth}
    \centering
    \includegraphics[width= 0.8\textwidth]{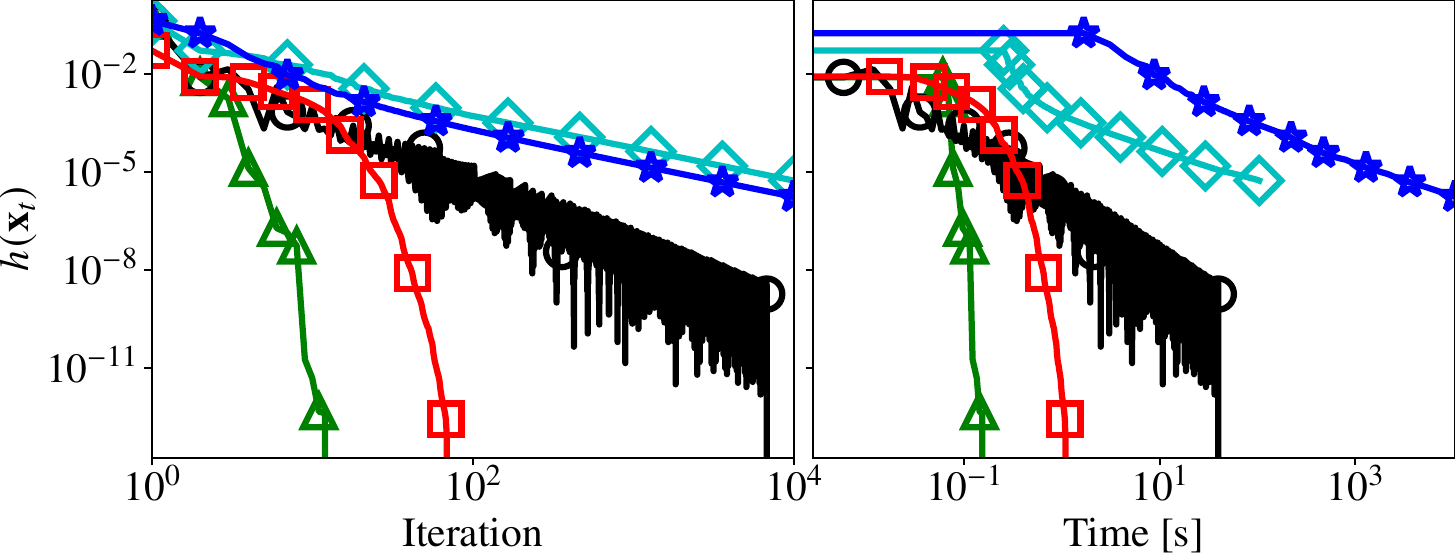}
    \includegraphics[width= 0.8\textwidth]{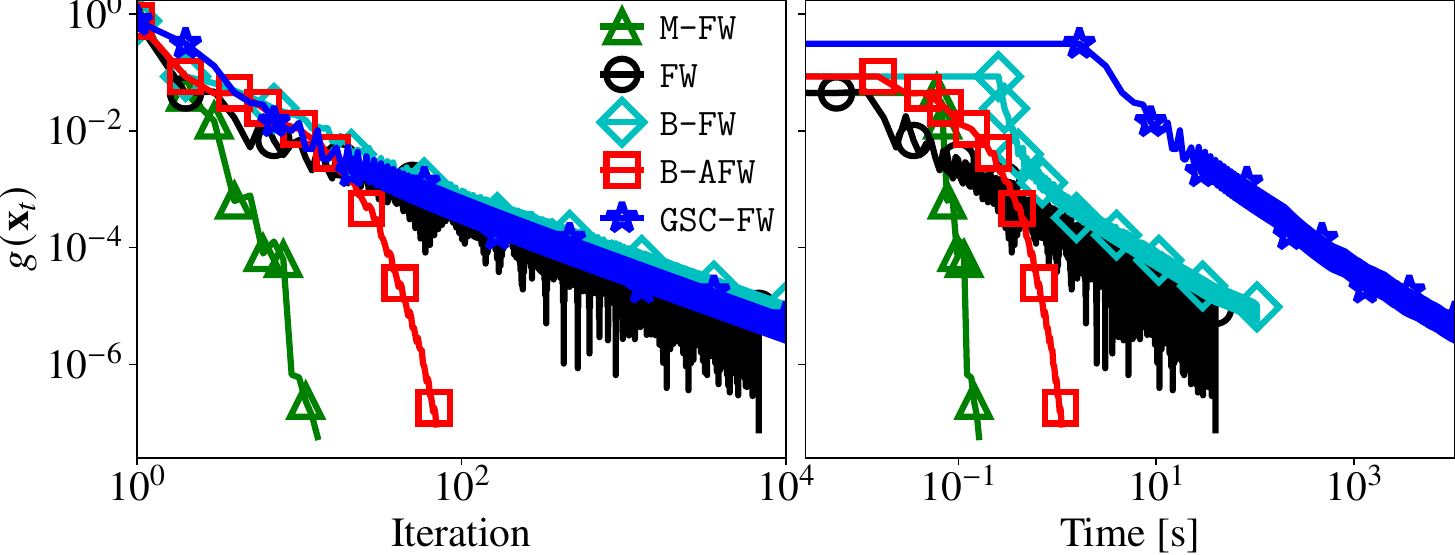}
      \caption{\texttt{a4a}: $(N,n) = (4781, 121)$}
\end{subfigure}
\hfill
\begin{subfigure}[b]{0.48\textwidth}
  \centering
  \includegraphics[width= 0.8\textwidth]{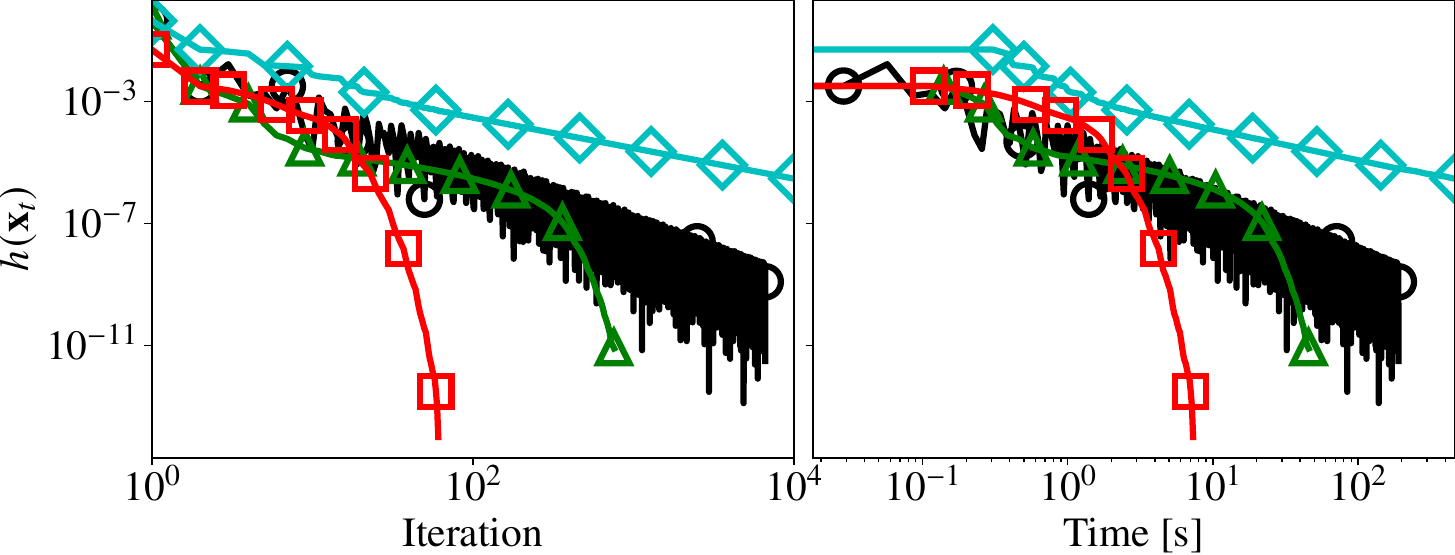}
  \includegraphics[width= 0.8\textwidth]{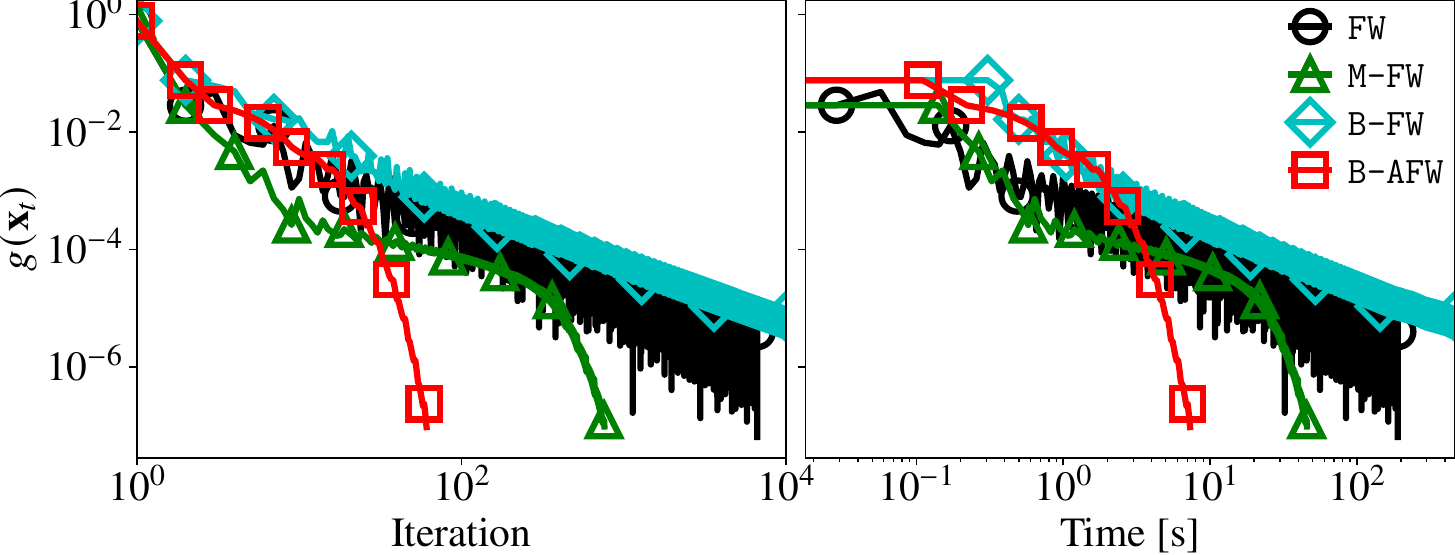}
  \caption{\texttt{a8a}: $(N,n) = (22696, 123)$}
\end{subfigure}
\caption{\textbf{Logistic regression}: Convergence of $h(\vx_t)$ and $g(\vx_t)$ vs. $t$ and wall-clock time for instances of the LIBSVM dataset.}
  \label{fig:logreggathered}
\end{figure}

\smallskip\noindent\textbf{Logistic regression.}
One of the motivating examples for the development of a theory of generalized self-concordant function is the logistic loss function, as it does not match the definition of a standard self-concordant function but shares many of its characteristics.
We consider a design matrix with rows $\va_i \in \mathbb{R}^n$ with $1 \leq i\leq N$ and a vector $\vy \in \{-1,1\}^N$
and formulate a logistic regression problem with elastic net regularization, in a similar fashion to \cite{liu2020newton}, with $f(\vx) = 1/N\sum_{i=1}^N \log(1 + \operatorname{exp}(-y_i \innp{\vx, \va_i})) + \mu/2 \norm{\vx}^2$, and $\cx$ is the $\ell_1$ ball of radius $\rho$.
The results can be seen in Figure~\ref{fig:logreggathered}.

\begin{figure}
  \centering
  \begin{tabular}{@{}c@{}}
    \includegraphics[width= 0.45\textwidth]{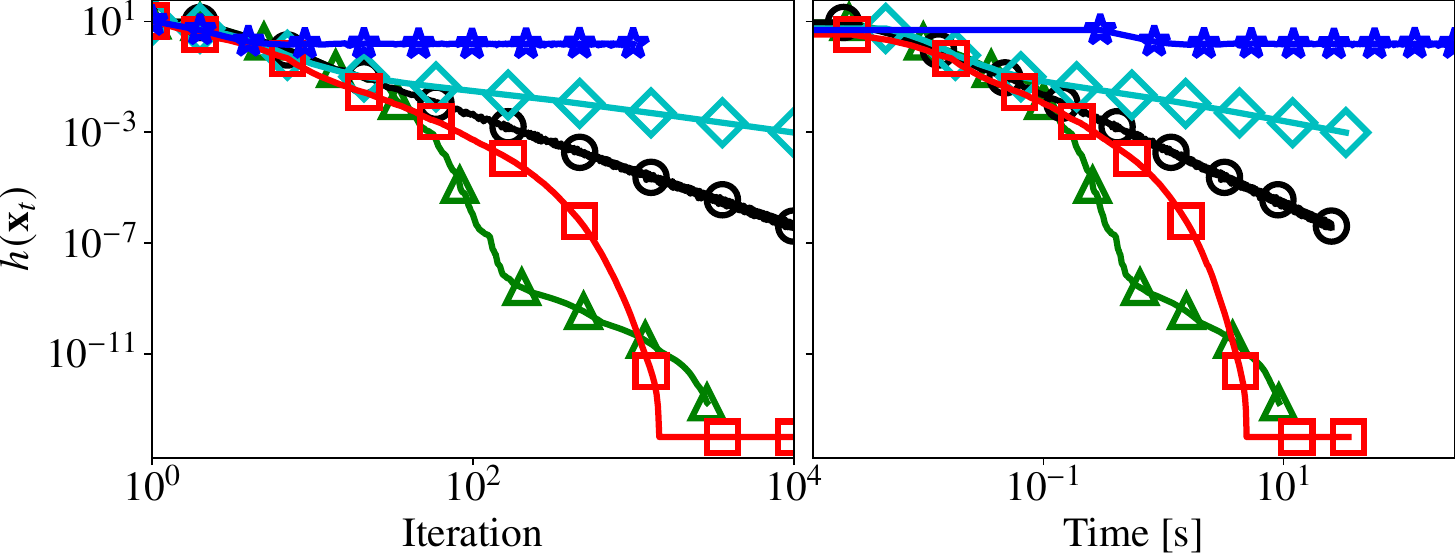} \\[\abovecaptionskip]
  \end{tabular}
  \begin{tabular}{@{}c@{}}
    \includegraphics[width= 0.45\textwidth]{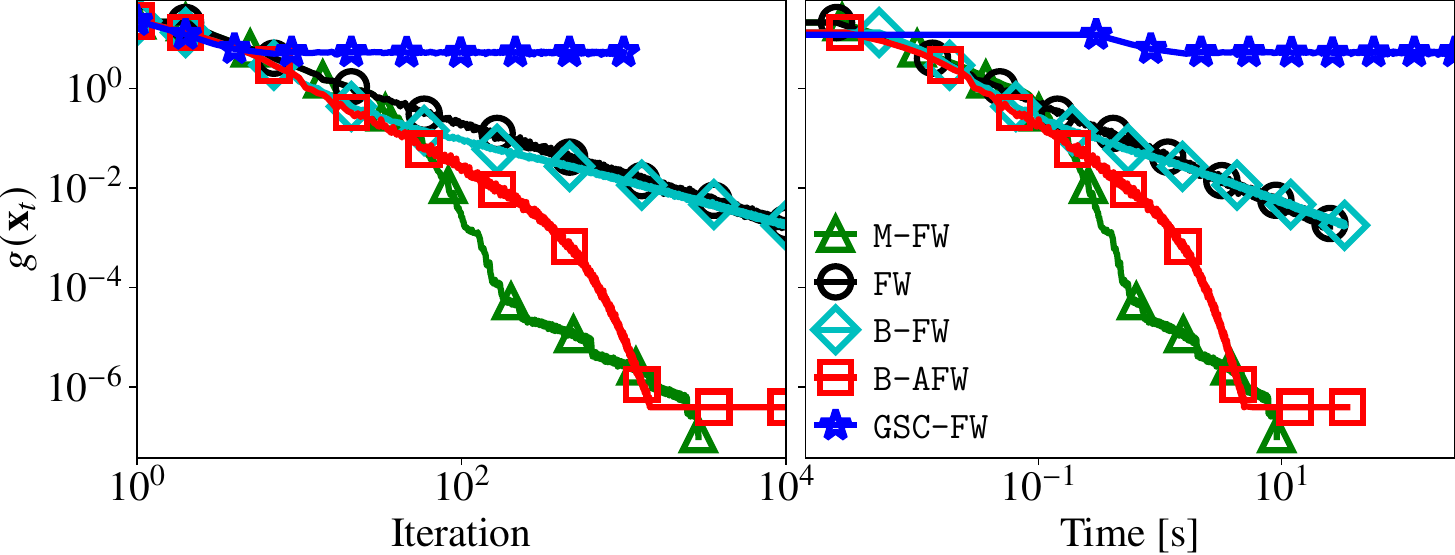} \\[\abovecaptionskip]
  \end{tabular}
  \vspace{-0.25cm}
  \caption{\textbf{Birkhoff polytope}: Convergence of $h(\vx_t)$ and $g(\vx_t)$ vs. $t$ and wall-clock time on \texttt{a2a}: $(N, n) = (2265, 114)$.}\label{fig:birkhoff}
\end{figure}

\smallskip\noindent\textbf{Birkhoff polytope.}
All applications previously considered all have in common a computationally inexpensive LMO that returns highly sparse vertices.
To complement the results, we consider a logistic regression problem over the Birkhoff polytope, where the LMO is implemented with the Hungarian algorithm, and is not as inexpensive as in the other examples.
We use a quadratic regularization parameter $\mu=100/\sqrt{N}$ where $N$ is the number of samples. The results are presented in Figure~\ref{fig:birkhoff}.

\paragraph{Monotonic step size: the numerical case}\label{par:numerics}
The computational experiments highlighted that the Monotonic Frank-Wolfe
performs well in terms of iteration count and time against other Frank-Wolfe and Away-step Frank-Wolfe variants.
Another advantage of a simple step size computation procedure is its numerical stability.
On some instances, an ill-conditioned gradient can lead to a plateau of the primal and/or dual progress.
Even worse, some step-size strategies do not guarantee monotonicity and can result in the primal value increasing over some iterations.
The numerical issue that causes this phenomenon is illustrated by running the methods of
the \texttt{FrankWolfe.jl} package over the same instance using $64$-bits floating-point numbers
and Julia \texttt{BigFloat} types (which support arithmetic in arbitrary precision to remove numerical issues).

\begin{figure}
  \centering
  \begin{subfigure}[b]{0.48\textwidth}
      \centering
      \includegraphics[width= 0.8\textwidth]{Images/Primal_gap_Portfolio_appx2_normal.pdf}
      \includegraphics[width= 0.8\textwidth]{Images/Dual_gap_Portfolio_appx2_normal.pdf}
        \caption{64-bit floating point}
  \end{subfigure}
  \hfill
  \begin{subfigure}[b]{0.48\textwidth}
    \centering
    \includegraphics[width= 0.8\textwidth]{Images/Primal_gap_Portfolio_appx1_big.pdf}
    \includegraphics[width= 0.8\textwidth]{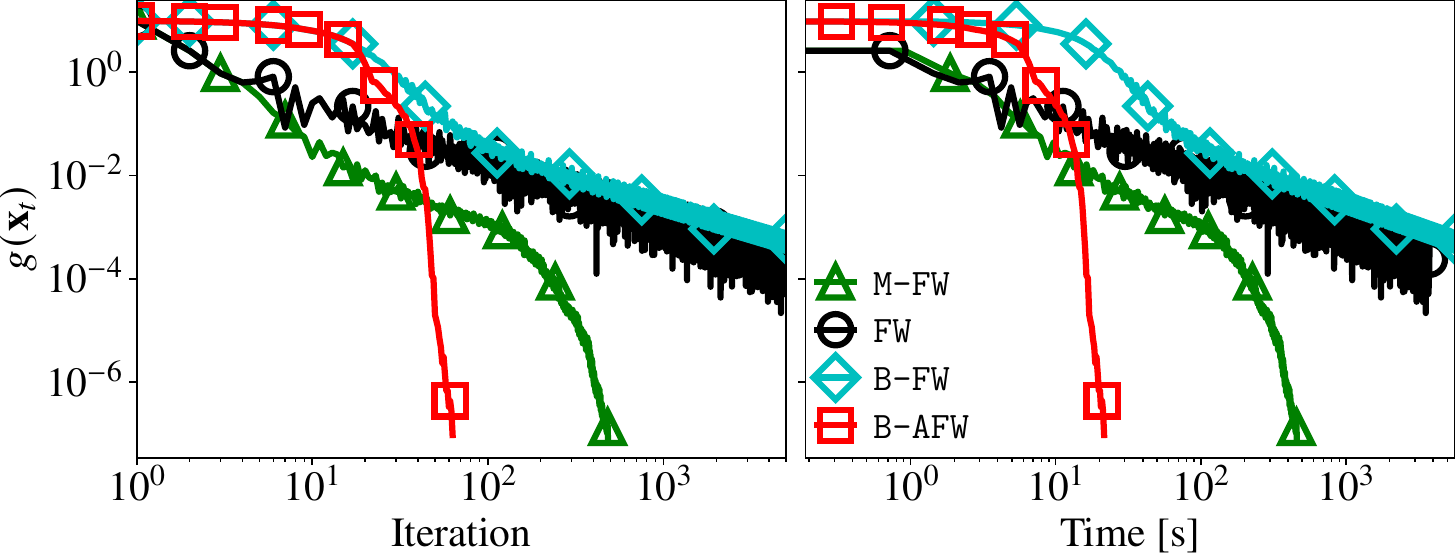}
    \caption{Arbitrary precision using \texttt{BigFloat}}
  \end{subfigure}
  \caption{Ill-conditioned portfolio optimization problem in dimension $n=2000$.}
    \label{fig:illcondition}
  \end{figure}



In \cref{fig:illcondition}, we compare the primal and dual gap progress of different algorithms on a portfolio instance.
In the finite-precision execution, we observe a plateau of the dual gap for both \texttt{M-FW} and \texttt{B-AFW}.
The primal value however worsens after the iteration where \texttt{B-AFW} reaches its dual gap plateau.
In contrast, \texttt{M-FW} reaches a plateau in both primal and dual gap at a certain iteration.
Note that the primal value at the point where the plateau is hit is already below $\sqrt{\varepsilon_{\text{float}64}}$,
the square root of the machine precision.
In arbitrary-precision arithmetic,
instead of reaching a plateau or deteriorating, \texttt{B-AFW} closes the dual gap tolerance and terminates before other methods.
Although this observation (made on several instances of the portfolio optimization problem) only impacts ill-conditioned problems,
it suggests \texttt{M-FW} may be a good candidate for a numerically robust default implementation of Frank-Wolfe algorithms.

\section{A stateless simple step variant}
The simple step-size strategy presented in \cref{alg:halvingmfw} ensures monotonicity and
domain-respecting iterates by maintaining a ``memory'' $\phi_t$ which is the number of performed halvings.
The number of halvings to reach an accepted step is bounded, but the corresponding factor $2^{\phi_t}$
is carried over in all following iterations, which may slow down progress.
We propose an alternative step size that still ensures the monotonicity and domain-preserving properties,
but does not carry over information from one iteration to the next.

\begin{algorithm}[H]
	\caption{Stateless Monotonic Frank-Wolfe}
  \label{alg:statelessfw}
  \footnotesize
\begin{algorithmic}[1]
\Require Point $\vx_0\in \cx \cap \dom{f}$, function \(f\)
\Ensure Iterates $\vx_1, \dotsc \in \cx$
\hrulealg
\For{$t=0$ \textbf{to} $\dotsc$}
  \State$\mathbf{v}_t\leftarrow\argmin_{\mathbf{v}\in\cx}\innp{\nabla f(\vx_t),\mathbf{v}}$
  \State$\gamma_t \leftarrow 2/(t+2)$
  \State$\vx_{t+1}\leftarrow \vx_t+\gamma_t(\mathbf{v}_t-\vx_t)$
  \While{$\vx_{t+1} \notin \dom f$ \textbf{ or } $f(\vx_{t+1}) > f(\vx_t)$}
      \State$\gamma_t \leftarrow \gamma_t /2$
      \State$\vx_{t+1}\leftarrow \vx_t+\gamma_t(\mathbf{v}_t-\vx_t)$
  \EndWhile
\EndFor
\end{algorithmic}
\end{algorithm}

Note that \cref{alg:statelessfw}, presented above, is stateless since it is equivalent to \cref{alg:halvingmfw} with $\phi_t$ reset to zero between every outer iteration.
This resetting step also implies that the per-iteration convergence rate of the stateless step is at least as good as the simple step,
at the potential cost of a bounded number of halvings, with associated ZOO and DOO calls at each iteration.
Finally, we point out that the stateless step-size strategy can be viewed as a particular instance of a backtracking line search
where the initial step size estimate is the agnostic step size $2/(t+2)$.

\begin{figure}
\centering
\begin{subfigure}[b]{0.48\textwidth}
    \centering
    \includegraphics[width= 0.8\textwidth]{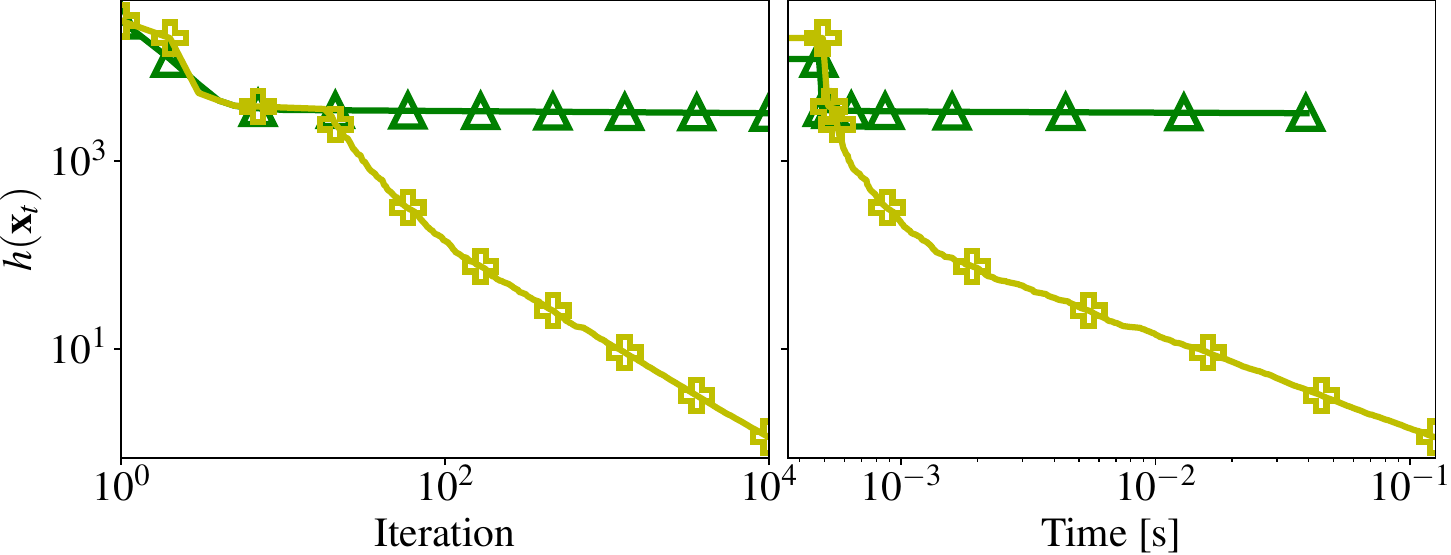}
    \includegraphics[width= 0.8\textwidth]{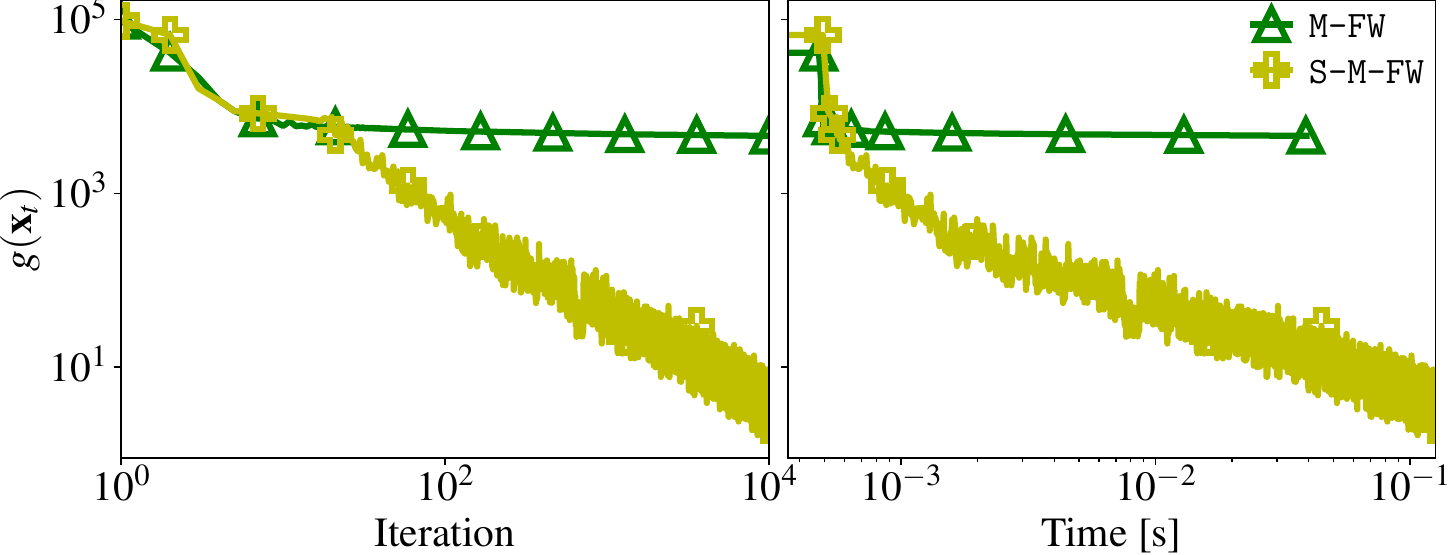}
      \caption{$n=30$}
\end{subfigure}
\hfill
\begin{subfigure}[b]{0.48\textwidth}
  \centering
  \includegraphics[width= 0.8\textwidth]{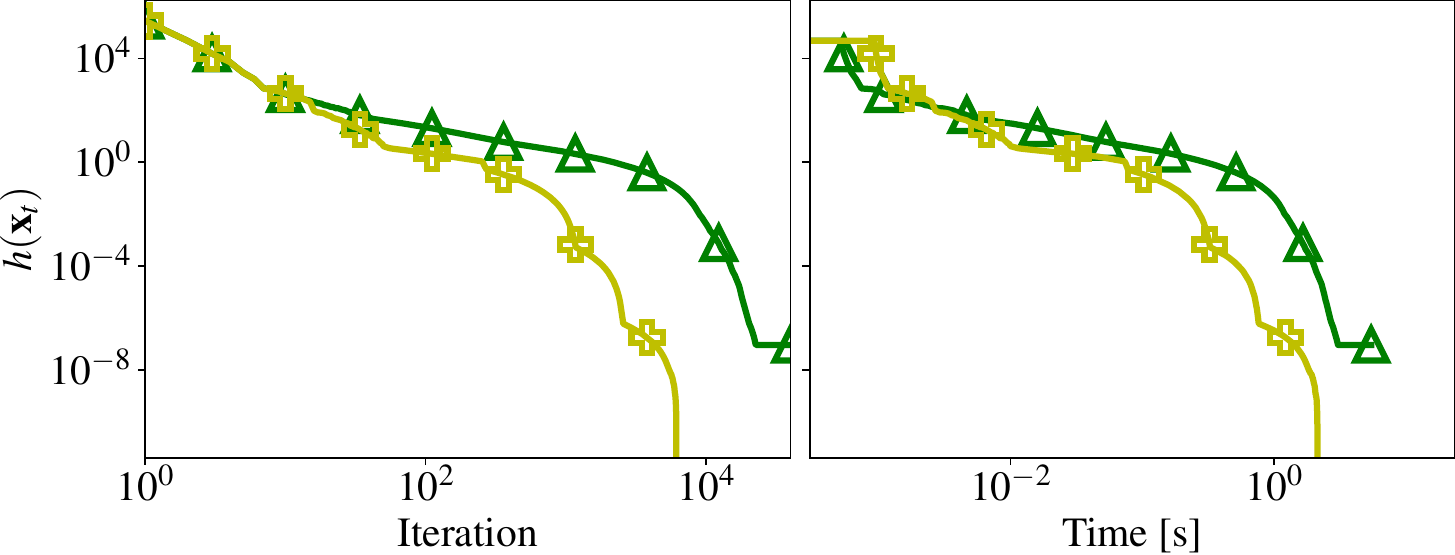}
  \includegraphics[width= 0.8\textwidth]{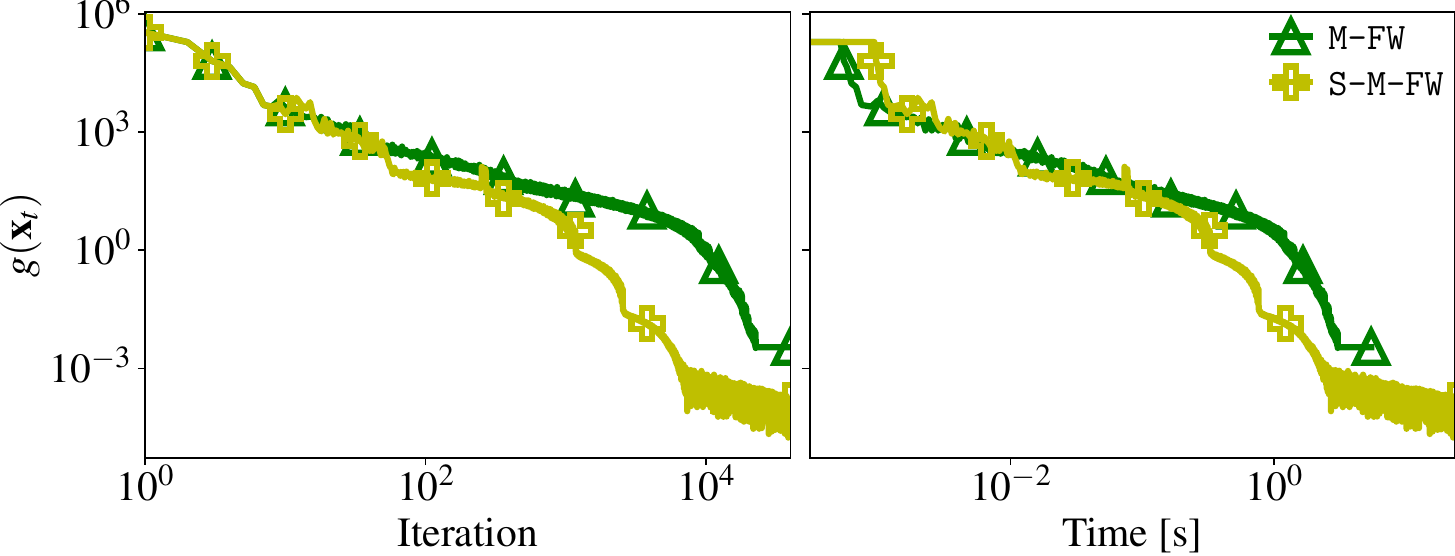}
  \caption{$n=1000$}
\end{subfigure}
\caption{\textbf{Stateless step size}: comparison of the stateless and simple steps on badly conditioned problems.}
  \label{fig:statelessgathered}
\end{figure}

We compare the two strategies on a random and badly conditioned problem with objective function:
\begin{equation*}
  x^\top Q x + \innp{b, x} + \mu \varphi_{\mathbb{R}_+}(x)
\end{equation*}
where $Q$ is a symmetric positive definite matrix with log-normally distributed eigenvalues and $\varphi_{\mathbb{R}_+}(\cdot)$
is the log-barrier function of the positive orthant. We optimize instances of this function over the $\ell_1$-norm ball
in dimension $30$ and $1000$.
The results are shown in Figure~\ref{fig:statelessgathered}. On both of these instances, the simple step progress is slowed down or even seems stalled in comparison to the stateless
version because a lot of halving steps were done in the early iterations for the simple step size, which penalizes progress over the whole run.
The stateless step-size does not suffer from this problem, however, because the halvings have to be performed at multiple iterations when using the stateless step-size strategy,
the per iteration cost of the stateless step-size is about three times that of the simple step-size.
Future work will consider additional restart conditions, not only on $\phi_t$ of \cref{alg:halvingmfw}, but also
on the base step-size strategy employed, similar to \citet{kerdreux2019restarting}.

\section*{Conclusion}

We introduced FQ variants based on open-loop
step sizes $\gamma_t = 2/(t+2)$ to obtain a $\mathcal{O}(1/t)$
convergence rate for generalized self-concordant functions in terms of primal and
FW gaps. This algorithm neither requires second-order information, nor line searches, and allows us to bound the number of zeroth-, first-order oracle, domain oracle, and linear oracle calls
needed to obtain the target accuracy.
We also show improved convergence rates for several variants in various cases of interest and prove that the AFW \citep{wolfe70,lacoste15} and BPCG \cite{tsuji2022pairwise} algorithms coupled with the backtracking line search of \cite{pedregosa2018step} can achieve linear convergence rates over polytopes when minimizing generalized self-concordant functions.

\section*{Acknowledgements}

Research reported in this paper was partially supported through the Research Campus Modal funded by the German Federal Ministry of Education and Research (fund numbers 05M14ZAM,05M20ZBM) and the Deutsche Forschungsgemeinschaft (DFG) through the DFG Cluster of Excellence MATH+. We would like to thank the anonymous reviewers for their suggestions and comments.


\bibliography{refs}
\bibliographystyle{icml2021}

\end{document}